%% file: main.tex
\documentclass[10pt]{article} 
\include{0Prae}

\title{Tensor approximation of the self-diffusion matrix of tagged particle processes}
\usepackage{authblk}
\author[1]{Jad Dabaghi}
\author[2,3]{Virginie Ehrlacher}
\author[4]{Christoph Strössner}
\affil[1]{{\footnotesize L\'eonard de Vinci P\^ole Universitaire, Research Center, 92916 Paris La D\'efense, France}}
\affil[2]{{\footnotesize Ecole des Ponts ParisTech, Marne-la-Vall\'ee, France}}
\affil[3]{{\footnotesize INRIA Paris, France, MATHERIALS team-project}}
\affil[4]{{\footnotesize \'Ecole Polytechnique F\'ed\'erale de Lausanne (EPFL), Institute of Mathematics, CH-1015 Lausanne, Switzerland}}
\date{}

\begin{document}
\maketitle

{\centering \vspace{-0.8cm} {\small
\url{jad.dabaghi@devinci.fr}, \url{virginie.ehrlacher@enpc.fr}, \url{christoph.stroessner@epfl.ch}} \\ \vspace{1cm}
\centering {\date{\today}} \\ \vspace{1cm}
}

\begin{abstract}
The objective of this paper is to investigate a new numerical method for the approximation of the self-diffusion matrix of a tagged particle process defined on a grid. While standard numerical methods make use of long-time averages of empirical means of deviations of some stochastic processes, and are thus subject to statistical noise, we propose here a tensor method in order to compute an approximation of the solution of a high-dimensional quadratic optimization problem, which enables to obtain a numerical approximation of the self-diffusion matrix. The tensor method we use here relies on an iterative scheme which builds low-rank approximations of the quantity of interest and on a carefully tuned variance reduction method so as to evaluate the various terms arising in the functional to minimize. In particular, we numerically observe here that it is much less subject to statistical noise than classical approaches.
\end{abstract}

\section{Introduction}

The aim of this paper is to propose a new numerical method for the approximation of the self-diffusion matrix of a tagged particle process on a grid~\cite{komorowski2012fluctuations} which is less subject to statistical noise than standard Monte-Carlo strategies. 

\medskip

The coefficients of the self-diffusion matrix of this process are defined as long-time averages of the mean-square displacement of the tagged particle~\cite{demasi2006mathematical,kipnis1994hydrodynamical,Quastel92,Blondel18}. The knowledge of the latter enables to identify the hydrodynamic limits of multi-species symmetric exclusion processes~\cite{Quastel92} and is thus of tremendous importance for applications. However, computing a numerical approximation of these coefficients is not an easy. \corr{Indeed, classical approaches consist in truncating the computational domain to a finite-size supercell with periodic boundary conditions and sampling a large number of realisations of the tagged particle trajectories in order to approximate the mean-square displacement by an empirical average computed with a standard Monte-Carlo approach. Moreover, the value of a finite final time has to be chosen beforehand so as to approximate the long-time limit. In~\cite{landim2002finite}, it has been proved that the error linked to the truncation of the computational domain decays exponentially with the size of the supercell. In contrast, the statistical error of the approximation linked to the use of a finite number of random samples of the trajectories of the tagged particle decays as the inverse of the square root of the number of samples. As a consequence, the main source of error in the practical computation of approximations of the self-diffusion matrix is due to statistical noise~\cite{Toninelli04,Thorneywork15,Hao09, Ferrando93,Sanz10,Deutsch91}.}

\medskip

In this paper, we propose an alternative numerical method based on the fact that the self-diffusion coefficients can be equivalently reformulated using the solution of a deterministic high-dimensional optimization problem~\cite{Blondel18,Quastel92, landim2002finite}. The main mathematical ingredient of the numerical method we investigate here is the use of low-rank tensor decompositions to obtain a numerical approximation of this solution, which enables to bypass the curse of dimensionality. This approach leads to an iterative scheme, each elementary step of which boils down to a  minimization problem over the set of low-rank tensors, which is solved using a classical alternating linear scheme~\cite{Beylkin05,Grasedyck13,Rohwedder13}.
Our numerical experiments demonstrate that this low-rank approach leads to very accurate approximations of the self-diffusion coefficients on finite-size grids.

As an illustration of the method, we use the obtained numerical approximations in order to compute the self-diffusion matrix of a two-dimensional tagged particle process defined on a Cartesian grid. \corr{We see this work as a preliminary step to computing diffusion coefficients out of Kinetic Monte Carlo
simulations for practical applications.}

\medskip

This article is organized as follows. In Section~\ref{sec:setting}, we introduce the lattice-based stochastic hopping model, its hydrodynamic limit and the two equivalent definitions of the self-diffusion matrix of the tagged particle process our work is based upon. The numerical approach we propose here in order to compute a numerical approximation  is presented in Section~\ref{sec:low:rank}. Finally, the efficiency of our approach is illustrated through several numerical experiments presented in~Section~\ref{sec:numerical:experiments}.

\section{Self-diffusion matrix}\label{sec:setting} 

\subsection{Infinite-dimensional definition} 

Let $d=1,2,3$ denote the dimension of the problem. We consider a symmetric tagged particle process defined on the infinite grid $\mathbb{Z}^d$. Let $K\in \mathbb{N}^*$ denote the number of possible jump directions for the particles and $(\bv_k)_{1\leq k \leq K} \subset \mathbb{Z}^d \setminus\{0\}$ the set of possible jump directions. For all $1\leq k \leq K$, the probability of jumping in the direction $\bv_k$ is denoted by $p_k \in \left]0,1\right]$. This jumping scheme is assumed to be symmetric in the sense that if the jump in the direction $\bm v_k$ occurs with probability $p_k$, then the jump in the direction $- \bm v_k$ occurs with the same probability.

The self-diffusion matrix application
$$
\matD_s : \left\{ 
\begin{array}{ccc}
[0,1] & \to & \mathbb{R}^{d\times d}\\
\rho & \mapsto & \matD_s(\rho):= \left(\matD_{s,ij}(\rho)\right)_{1\leq i,j \leq d}\\
\end{array}\right.
$$
can be defined in the following two ways.  

\medskip

\paragraph{Definition as optimization problem.}
Let us first introduce some notation. Let
$S:= \mathbb{Z}^d \setminus \{\bm 0\}$. For all $\bEta:= (\eta_\bs)_{\bs\in S} \in \{0,1\}^S$ and all $\by\neq \bz \in S$, we define by $\bEta^{\by, \bz}:= (\eta^{\by, \bz}_\bs)_{\bs\in S}$ the element of $\{0,1\}^S$ such that
$$
\eta^{\by, \bz}_\bs:=\left\{
\begin{array}{ll}
\eta_\bs & \mbox{if } \bs\neq \by, \bz,\\
\eta_\by& \mbox{ if }\bs = \bz,\\
\eta_\bz & \mbox{ if } \bs = \by.\\
\end{array}\right. 
$$
Furthermore, for all $\bw \in S$, we define by $\bEta^{\bm 0, \bw}:= (\eta^{\bm 0, \bw}_\bs)_{\bs\in S}$ the element of $\{0,1\}^S$ such that
$$
\eta^{\bm 0, \bw}_\bs:=\left\{
\begin{array}{ll}
\eta_{\bs+\bw} & \mbox{if } \bs\neq -\bw, \\
0& \mbox{ if }\bs = -\bw.\\
\end{array}\right. 
$$For \corr{a mean particle density} $\rho \in [0,1]$, we denote by $\rho^\otimes$ the Bernoulli product measure on $\{0,1\}^S$ with marginals given by
$$
\rho^\otimes\left(\left\{\bEta:=(\eta_{\bs})_{\bs\in S}: \eta_{{\hat{\bs}}} = 1\right\}\right) = \rho, 
$$
for all ${\hat{\bs}} \in S$. 
Let us also define the set of functions $H_\rho:= L^2_{\rho^\otimes}\left(\{0,1\}^S\right)$, where $L^2_{\rho^\otimes}$ denotes the $L^2$ Lebesgue space with measure $\rho^\otimes$.

For a \corr{drift direction} $\bu \in \R^d$ and \corr{mean particle density} $\rho\in [0,1]$, the self-diffusion coefficient $\bu^T\mathbb{D}_s\left(\rho\right)\bu$ is given by: 
\begin{equation}\label{eq:infinite}
\bu^T \matD_s(\rho) \bu:= {2} \mathop{\inf}_{\Psi\in H_\rho} \mathbb{E}_{\rho^\otimes}\left[\sum_{k=1}^K p_k
\bra{ (1- \eta_{\bv_k}) \bra{
\bu\cdot \bv_k + \Psi(\bEta^{\bm 0,\bv_k}) - \Psi(\bEta)}^2 + {\frac{1}{2}}
\sum_{\substack{\by \in S  \\ \by+\bv_k \neq \bm 0 } } \bra{
\Psi(\bEta^{\by+\bv_k,\by})-\Psi(\bEta)}^2 }  \right], 
\end{equation}
where the notation $\mathbb{E}_{\rho^\otimes}$ refers to the fact that the expectation is taken over the product measure $\rho^\otimes$~\cite{Quastel92,komorowski2012fluctuations}. Problem (\ref{eq:infinite}) thus reads as an infinite-dimensional optimization problem over the set $H_\rho$. 

\begin{rmrk} 
Naturally, for all $\rho\in [0,1]$, since $\matD_s(\rho)$ is a symmetric matrix, one can easily deduce the full matrix $\matD_s(\rho)$ from the knowledge of $\bu^T \matD_s(\rho)\bu$ for a few vectors $\bu \in \mathbb{R}^d$. 
\end{rmrk}

\medskip

\paragraph{Definition as long time mean square deviation.}
The quantity $\bu^T \matD_s(\rho) \bu$ can be equivalently expressed as the long time limit of the following expectation~\cite[Theorem 2.3]{Quastel92}. Let us assume that, at time $t=0$, the tagged particle is located at position $\bm 0$ and all other sites are occupied with probability $\rho$ following a Bernoulli distribution. The duration between two consecutive jumping events follow an exponential law with parameter $1$. The jumping directions (respectively rates) are given by $\{\bm v_1 , \dots , \bm v_K \}$ (respectively $p_1,\dots,p_K$). Jumps are not allowed if the final site of the jumping particle is already occupied. Let $\bm w(t)$ denote the position of the tagged particle at time $t\geq 0$. It then holds that, for a \corr{drift direction} $\bu \in \mathbb{R}^d$ and \corr{mean particle density }$\rho \in [0,1]$, the quantity $\bu^T \matD_s(\rho) \bu$ can be equivalently formulated as \corr{the long-time limit mean square deviation of the tagged particle in the direction $\bu \in \mathbb{R}^d$, i.e.}
\begin{equation}\label{eq:infiniteLimit}
\bu^T \matD_s(\rho) \bu = {2} \lim_{t \to \infty} \frac{\mathbb{E}_{\rho^\otimes} [ \langle \bu, \bm w(t)\rangle^2 ]}{t},
\end{equation}
where $\langle \cdot, \cdot \rangle$ denotes the euclidean scalar product of $\mathbb{R}^d$.

\begin{rmrk}
Starting from a Bernoulli-product measure and symmetric transition rates, it is a classical problem in probability theory to study the motion of a tagged particle on $\mathbb Z^d$~\cite{Saada87,Kipnis86,Liggett99}, see also~\cite{Chen19b,Gantert20} for recent work on tagged particles on trees. For $d \geq 2$ and symmetric nearest neighbor transition rates, the tagged particle is known to satisfy a central limit theorem with non-degenerate limiting variance. However, to our best knowledge, a general closed formula for the limiting variance is not available. 

Notice that when starting from a Bernoulli-product measure with a tagged particle in the origin, the resulting environment process is stationary with respect to the Bernoulli-product measure conditioned to contain a particle in the origin, usually called the Palm measure. This suggests that the limiting variance for the tagged particle can be described in terms of the sum of separable functions, i.e. in terms of a low-rank function. In turn, this indicates that the equivalent characterization of the limiting variance as in Equation~\eqref{eq:infinite} is also related to low-rank functions. \end{rmrk}

\medskip

Let us point out that the quantity $\bu^T \matD_s(\rho) \bu$ cannot be computed exactly in practice, neither using expression (\ref{eq:infinite}) nor expression (\ref{eq:infiniteLimit}). On the one hand, (\ref{eq:infinite}) reads as an infinite-dimensional optimization problem and has to be approximated by a finite-dimensional optimization problem in practice~\cite{landim2002finite,Jara06}. On the other hand, (\ref{eq:infiniteLimit}) requires the computation of the long-time average of the stochastic process defined on the infinite lattice grid $\mathbb{Z}^d$. Equivalently, it has to be approximated in a finite-dimensional setting, as long-time average of the mean-square deviation of the tagged particle associated to a stochastic process defined on a finite-size grid with periodic boundary conditions. Both finite-dimensional approximations are presented in detail in the next section.

\subsection{Finite-dimensional approximation}

\begin{figure}
    \centering
    \includegraphics[width=0.8\textwidth]{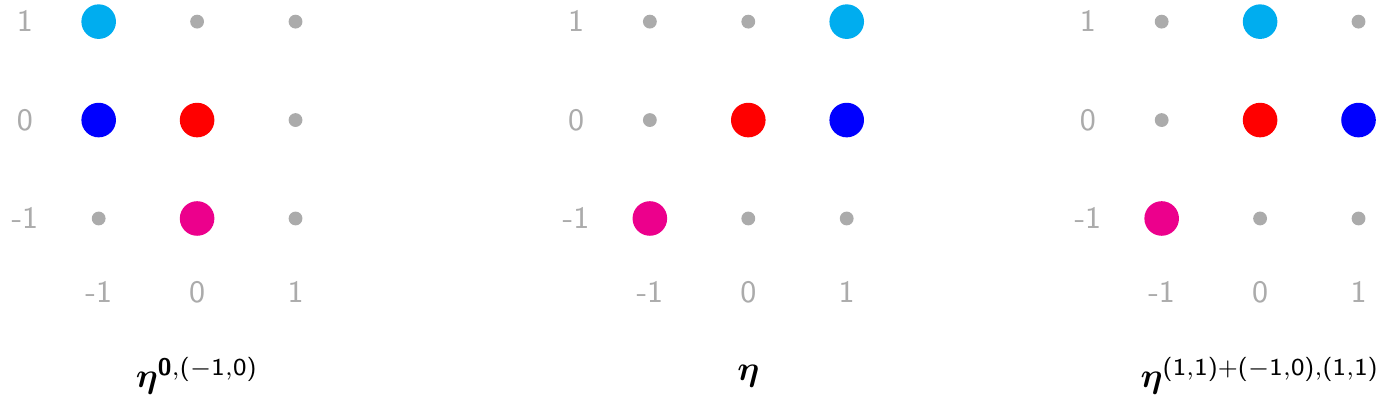}
    \caption{
    Middle: Visualization of one particular $\bEta \in S_1$ for $d=2$ with three occupied sites marked in blue, cyan and magenta. Additionally, we mark the tagged particle $\bm 0$ in red.
    Left: Visualization of the occupied sites of $\bEta^{\bm 0, (-1,0)}$. This can be seen as a jump of the imaginary red particle one step to the left, followed by immediately relabeling of the sites such that red particle remains at $\bm 0$. By exploiting the periodicity, we obtain $\bEta^{\bm 0, (-1,0)} \in S_1$.
    Right: Visualization of the occupied sites of $\bEta^{(1,1)+(-1,0),(1,1)}$. This can be seen as jump of the cyan particle one step to the left.
    }
    \label{fig:Jumps}
\end{figure}

\paragraph{Discretized minimization problems.}
Let $M \in \mathbb{N}^*$  denote a discretization parameter and introduce the finite grid $S_M:= \{-M, \cdots, M\}^d \setminus \{\bm 0\}$. For the sake of simplicity, we assume that $M \geq \max_{k \in \{1,\dots,K\}} ||v_k||_1$, where $||\cdot ||_1$ denotes the taxicab norm. For any $\bEta \in \{0,1\}^{S_M}$, we can construct by periodicity an extension $\widetilde{\bEta}:= (\widetilde{\eta}_\bs)_{\bs\in S} \in \{0,1\}^S$ by assuming with a slight abuse of notation that the site $\bm 0$ is occupied, i.e. $\widetilde{\eta}_\bs = 1$ for $\bs \in (2M+1)\cdot\mathbb{Z}^d \setminus \{\bm 0\}$. Using this notation, for all $\bEta\in \{0,1\}^{S_M}$ and all $\by, \bz, \bw \in S$, we define $\bEta^{\by, \bz}\in \{0,1\}^{S_M}$ and $\bEta^{\bm 0, \bw}\in \{0,1\}^{S_M}$ as 
$$
\bEta^{\by, \bz}:= \left( \widetilde{\eta}^{\by, \bz}_\bs \right)_{\bs \in S_M} \quad \mbox{ and }\quad \bEta^{\bm 0, \bw}:= \left( \widetilde{\eta}^{\bm 0, \bw}_\bs \right)_{\bs \in S_M}.
$$
Figure~\ref{fig:Jumps} shows an illustration of $\bEta^{\by, \bz}$ and $\bEta^{\bm 0,\bw}$.

Let $N:= (2M+1)^d -1 = {\rm Card}(S_M)$. 
For all $\ell \in \{0,\dots,N\}$, let $C_{M,\ell} := \{ \bEta \in \{0,1\}^{S_M} | \sum_{\bs \in S_M} \eta_\bs = \ell\}$ denote the set of all possible configurations of the particles on $S_M$ so that the total number of occupied sites is equal to $\ell$. 

Let us define $H_M:= \left\{ \Psi: \{0,1\}^{S_M} \to \mathbb{R}\right\}$. For every $\bu \in \R^d$ and $l \in \{0,\dots,N\}$, we introduce  the quadratic functional $A_{M,\ell}^\bu : H_M \to \mathbb{R}$ defined by
\begin{equation*}
A_{M,\ell}^{\bm u}(\Psi) := \frac{1}{|C_{M,\ell}|} \sum_{\bEta \in C_{M,\ell}} 
\sum_{k=1}^K p_k
\bra{ (1- \eta_{\bv_k}) \bra{
\bu\cdot \bv_k + \Psi(\bEta^{\bm 0,\bv_k}) - \Psi(\bEta)}^2 + {\frac{1}{2}}
\sum_{\substack{\by \in S_M \\ \by+\bv_k \neq \bm 0 } } \bra{
\Psi(\bEta^{\by+\bv_k,\by})-\Psi(\bEta)}^2   }.
\end{equation*}
    
Then, assuming that $\overline{\rho} = \frac{\ell}{N}$ for some $0\leq \ell \leq N$, one can define~\cite{Blondel18} for all $\bu \in \mathbb{R}^d$, 
\begin{equation}
\label{eq:approx:Ds}
\bu^T \mathbb{D}^M_s\left(\frac{\ell}{N}\right) \bu:= {2} \umin{\Psi \in H_M} A_{M, \ell}^{\bm u} (\Psi).
\end{equation}

It is proved in~\cite{landim2002finite} that $\displaystyle \mathop{\lim}_{\begin{array}{c}
    M\to +\infty \\
    \frac{\ell}{N} \to \overline{\rho} \\
    \end{array}}
\bu^T \mathbb{D}^M_s\left(\frac{\ell}{N}\right) \bu = \bu^T \mathbb{D}_s(\overline{\rho}) \bu$. 

\medskip

\paragraph{Combined minimization problem.}
The collection of sets $C_{M,0},\dots,C_{M,N}$ form a partition of the set $\{0,1\}^{S_M}$. Observe for a given $\Psi\in H_M$, $A_{M,\ell}^{\bm u}(\Psi)$ only depends on the values of $\Psi(\bEta)$ for $\bEta \in C_{M,\ell}$, since $\bEta \in C_{M,\ell}$ implies that also $\bEta^{\bm 0,\bm v_k} \in C_{M,\ell}$ and $\bEta^{\bm y + \bm v_k, \bm y} \in C_{M,\ell}$ for all $1\leq k \leq K$.
As a consequence,  if $\Psi^{M,\bm u}_{\textsf{opt}}\in H_M$ is a minimizer of 
\begin{equation}\label{eq:minopt}
\min_{\Psi \in H_M} A^{\bm u}_M (\Psi),
\end{equation}
where 
\begin{equation}
    \label{eq:def:Au}
A_M^{\bm u}(\Psi) := \sum_{\bEta \in \{0,1\}^{S_M}} 
\sum_{k=1}^K p_k
\bra{ (1- \eta_{\bv_k}) \bra{
\bu\cdot \bv_k + \Psi(\bEta^{\bm 0,\bv_k}) - \Psi(\bEta)}^2 + {\frac{1}{2}}
\sum_{\substack{\by \in S_M \\ \by+\bv_k \neq \bm 0 } } \bra{
\Psi(\bEta^{\by+\bv_k,\by})-\Psi(\bEta)}^2 },
\end{equation}
then it holds that $A_{M,\ell}^{\bm u}(\Psi^{M,\bm u}_{\textsf{opt}}) = \min_{\Psi \in H_M} A_{M,\ell}^{\bm u}(\Psi)$ for all $\ell \in \{0,\dots,N\}$. The knowledge of $\Psi^{M,\bm u}_{\textsf{opt}}$ then allows us to compute $\bu^T \mathbb{D}_s^M\left(\frac{\ell}{N}\right) \bu$ for all $0\leq \ell \leq N$ as $2 A_{M,\ell}^{\bm u} (\Psi^{M,\bm u}_{\textsf{opt}})$. Note that the minimization problem (\ref{eq:minopt}) is then independent of $\ell$, in contrast to (\ref{eq:approx:Ds}). 
However, the minimization problem~\eqref{eq:minopt} has $2^N$ degrees of freedom and is thus intractable for large values of $N$. 

The main contribution of our paper is to propose a low-rank tensor approximation algorithm to compute an approximation of $\Psi_{\rm opt}^{M,u}$ in order to mitigate this issue. Our proposed method is presented in Section~\ref{sec:low:rank}.

\medskip

\paragraph{Estimation of long-time mean square deviation.}
The quantity $\bu^T \mathbb{D}^M_s\left(\frac{\ell}{N}\right) \bu$ can be equivalently expressed as the long-time mean square deviation of a tagged particle evolving in a periodic environment~\cite{landim2002finite}. The Monte-Carlo algorithm is the standard method of choice in practice to compute an approximation of $\bu^T \mathbb{D}^M_s\left(\frac{\ell}{N}\right) \bu$ by means of empirical averages. The quality of the obtained approximations then depends on the choice of two numerical parameters, namely the number of Monte-Carlo samples (i.e. stochastic realizations of the periodized tagged particle process) and the value of the chosen final time.  

We aim to compare the low-rank tensor approximation algorithm we propose in this work with such Monte-Carlo methods. To this aim, we detail below the algorithm used to compute such empirical averages, and which we will use as a comparison with the method we develop here and which is presented in Section~\ref{sec:low:rank}.

Let $1 \leq \ell \in \leq N$. 
An initial environment $\bEta$ is obtained by sampling uniformly from $C_{M,\ell}$. We sample from exponential distributions to determine the next time a particle in the environment or the tagged particle performs a jump. 
Whenever a jump is performed the environment $\bEta$ is updated accordingly. 
Let $\bm w = \bm 0 \in \mathbb Z^d$ denote the initial position of the tagged particle. Throughout the simulation, we track the position of the tagged particle in $\mathbb Z^d$, i.e. whenever the tagged particle successfully jumps in direction $\bm v_k$, we set $\bm w = \bm  w +  \bm v_k$. 
To approximate the expectation of $\bm w$, we repeat this simulation $N_s$ times with different samples. Each of these simulations is stopped at the same stopping time $T$ at which we compute the approximation of $\bu^T \matD_s({\rho}) \bu$. This approach is formalized in Algorithm~\ref{alg:StochasticAlgorithm}.

\begin{algorithm}
\caption{Long-term Monte Carlo estimation}\label{alg:StochasticAlgorithm}
\begin{algorithmic}[1]
\algdef{SE}[SUBALG]{Indent}{EndIndent}{}{\algorithmicend\ }%
\algtext*{Indent}
\algtext*{EndIndent}
\State \textbf{Input}: $M\in \mathbb{N}^*$, $\bm u \in \mathbb{R}^d$, number of occupied sites $1\leq \ell \leq N$ with $N = (2M+1)^d -1$, final time $T>0$,  $\hat{N}_s \in \mathbb{N}^*$
\State \textbf{Output}: approximation $\alpha \approx \frac{1}{2} \bu^T \mathbb{D}^M_s\left(\frac{\ell}{N}\right) \bu = A_{M,\ell}^{\bm u} (\Psi^{M,\bm u}_{\textsf{opt}})$
\State $\bm w = \bm 0$, $\alpha = 0$, $N_s = \lceil \hat{N}_s/(\ell + 1)\rceil$
\For $i = 1,2,3,\dots,N_s$
    \State randomly initialize $\bEta \in C_{M,\ell}$ and set $t_{\textsc{total}} = 0$, 
    \While $\text{true}$
        \State sample $t_{\textsc{new}}$ from an exponential distribution with mean $\frac{1}{\ell + 1}$ and set $ t_{\textsc{total}}  = t_{\textsc{total}} + t_{\textsc{new}}$
        \If   $t_{\textsc{total}} > T$ 
        \State break
        \EndIf
        \State uniformly select either one occupied site $\bm y \in S_M$ with $\bEta_{\bm y} = 1$ or the tagged particle $\bm y = \bm 0$
        \State select a jump direction $\bm v_k$ with probability $p_k$
        \If $\tilde{\bEta}_{\bm y + \bm v_k} = \bm 0$, i.e. the target location is not occupied
            \State when $\bm y \neq \bm 0$ $\Rightarrow $ set $\bEta = \bEta^{\bm y + \bm v_k, \bm y}$ 
            \State when $\bm y = \bm 0$ $\Rightarrow $ set $\bEta = \bEta^{\bm 0, \bm v_k}$ and update $\bm w = \bm w + \bm v_k$ 
            \EndIf
    \EndWhile
    \State update $\alpha = \alpha + \langle \bm u, \bm w \rangle^2$ \label{line:UpdateAlphaBasedonU}
\EndFor
\State $\alpha = \frac{\alpha}{T N_s} $
\label{line:UpdateAlpha2}
\end{algorithmic}
\end{algorithm} 

\begin{rmrk}
In Algorithm~\ref{alg:StochasticAlgorithm}, we set $N_s = \lceil \hat{N}_s/(\ell + 1)\rceil$, where  $\hat{N}_s$ is a given input parameter. This ensures that the expected runtime is approximately equal for all choices of $\ell$. Additionally Figure~\ref{fig:SamplingStochasicVariance} shows that the variance of the output is comparable independently of $\ell$. 
If we were to use the same $N_s$ for all $\ell$, we would observe a much larger variance for smaller $\ell$ and the runtime would increase for larger $\ell$.
\end{rmrk}

\section{Low-rank solutions for the optimization problem}\label{sec:low:rank}

The aim of this section is to present the numerical method based on low-rank tensor approximations we propose for the resolution of the high-dimensional optimization problem (\ref{eq:minopt}). 
We first define low-rank functions in $H_M$ and introduce a fast and stable algorithm for the evaluation of $A_M^{\bm u}$. 
We then develop a successive minimization scheme to compute low-rank solutions of~\eqref{eq:minopt}. 
Each minimization step is performed using an alternating linear scheme, which relies on the fast and stable evaluations of $A_M^{\bm u}$.
Lastly, we discuss the evaluation of $A_{M,\ell}^{\bm u}$ for the computation of the self-diffusion coefficient~\eqref{eq:approx:Ds}.

\subsection{Low-rank functions}
A function $R \in H_M$ is called {a rank-$1$ or pure tensor product function} when it can be written as \[
R (\bEta) = \Pi_{\bs\in S_M} R_\bs(\eta_\bs), \quad {\forall \bEta = (\eta_\bs)_{\bs\in S_M}\in \{0,1\}^{S_M},}
\]
for some $R_\bs: \{0,1\} \to \R$ for $\bs\in S_M$. Let $H_M^1 \subset H_M$ denote the set of pure tensor product functions of $H_M$. 
Let $r \in \mathbb{N}$. A function $\Phi \in H_M$ is called a rank-$r$ function when it can be written as $\Phi = R^{(1)} + \dots + R^{(r)}$, where $R^{(k)} \in H_M^1$ for $1 \leq k \leq r$. We denote by $H_M^r \subset H_M$ the set of all rank-$r$ functions.
For all $r \in \mathbb{N}^*$, it holds 
\begin{equation*}
    \umin{\Psi \in H_M} A^{\bm u}_M(\Psi) = \umin{\Phi \in H^{2^N}_M} A^{\bm u}_M(\Phi) \leq \umin{\Phi \in H^{r+1}_M} A^{\bm u}_M(\Phi) \leq \umin{\Phi \in H^r_M} A^{\bm u}_M(\Phi) \leq \umin{R \in H^1_M} A^{\bm u}_M(R).
\end{equation*}
In the next sections, we derive an algorithm to compute an approximation of $\min_{\Phi \in H_M^r} A^{\bm u}_M(\Phi)$. 

\subsection{Fast and stable evaluation for low-rank functions}\label{sec:StableEvaluation}
In this section, we introduce a fast and stable method to evaluate $A^{\bm u}_M(\Psi)$ for $\Psi \in H^r_M$. 
A naive evaluation would require to sum over $2^N$ terms, which is intractable for large values of $N$.
In principle, the order of summation can be exchanged leading to terms of the form $\sum_{\bEta \in {\{0,1\}^{S_M}}} \Psi(\bEta)^2$, which can be evaluated individually.
However, subtracting such terms leads to a lot of numerical cancellation for $M > 1$. 
We circumvent these issues by treating the evaluation of $A^{\bm u}_M(\Psi)$ as the computation of the Frobenius norms of certain tensors. These Frobenius norms can then be evaluated in a fast and stable way.

\subsubsection{Reformulation as tensor norm }
For $\Psi \in H$, we consider the tensor $\mathcal{T}^{(\Psi)} \in \R^{2^N}$ with entries $\mathcal T^{(\Psi)}_{i_1,\dots,i_N} = \Psi(\bEta)$, where $\bEta = (\eta_{\bm s_1}, \dots, \eta_{\bm s_N})$ for some enumeration $\bm s_1, \dots, \bm s_N$ of the sites in $S_M$ and $\eta_{\bm s_j} = i_j-1$ for $1\leq j \leq N$. 
Note that $\mathcal{T}^{(\Psi)}$ is a tensor of at most rank-$r$ when $\Psi \in H^r_M$. In particular, when $\Psi$ is given in term of the functions $R_{\bm s}^{(k)}$ for $1 \leq j \leq N$, $1 \leq k \leq r$ , we can write $\mathcal{T}^{(\Psi)}$ in the so called CP-format~\cite{Kiers00,Kolda09} as 
\begin{equation}
    \label{eq:CPformat}
\mathcal T^{(\Psi)}_{i_1,\dots,i_N} = \sum_{k=1}^r \prod_{j=1}^N \bm a_{i_j}^{(j,k)},
\end{equation}
where $\bm a^{(j,k)} \in \R^2$  is defined as $\bm a^{(j,k)} = (R_{\bm s_j}^{(k)} (0),R_{\bm s_j}^{(k)}(1) )$. 

For $\bv,\by \in S_m$, we analogously define the tensor $\mathcal{T}^{(\Psi^{(\bm 0,\bm v)})} \in \R^{2^N} $ with entries $\mathcal T^{(\Psi^{(\bm 0,\bm v)})}_{i_1,\dots,i_N} = \Psi(\bEta^{\bm 0, \bv})$. When $\by + \bv \neq \bm 0$ we additionally define the tensor $\mathcal{T}^{(\Psi^{(\by + \bv, \by)})} \in \R^{2^N}$  with entries $\mathcal T^{(\Psi^{(\by + \bv, \by)})}_{i_1,\dots,i_N} = \Psi(\bEta^{\by + \bv, \by})$. We observe that these are again are rank-$r$ tensors when $\Psi \in H_M^r$. Lastly, we define the rank-$1$ tensor $\mathcal{T}^{(\bu \cdot \bv)} \in \R^{2^N}$ with entries $\mathcal{T}^{(\bu \cdot \bv)}_{i_1,\dots,i_N} = \bu \cdot \bv$.

For a tensor $\mathcal{T} \in \R^{2^N}$ and a site $\bm s$ we define the projection operators $\mathcal{P}_{\bm s}: \R^{2^N} \to \R^{2^N}$ as
\[
(\mathcal{P}_{\bm s} (\mathcal{T}))_{i_1,\dots,i_N} := \left\{
\begin{array}{ll}
\mathcal{T}_{i_1,\dots,i_N} & i_{\bm s} = 1  \\
0 & \, \textrm{otherwise} \\
\end{array}
\right., \quad \text{for } 1\leq i_1,\dots,i_N \leq 2,
\]
where $i_{\bm s}$ denotes to the index assigned to site $\bm s$. 
We can now rewrite~\eqref{eq:def:Au} as
\begin{equation}
    \label{eq:tensorformulation}
    A_M^{\bm u}(\Psi) = \sum_{k=1}^{K} p_k \Bigg (
    \norm{P_{\bm v_k}\big{(}\mathcal{T}^{(\bu \cdot \bv_k)} + \mathcal{T}^{(\Psi^{(\bm 0, \bv_k)})} - \mathcal{T}^{(\Psi^{(\by+\bv_k,\by)})}\big{)}}_F^2 + \frac{1}{2}
    \sum_{\substack{\by \in S_M \\ \by+\bv_k \neq \bm 0 } } \norm{  \mathcal{T}^{(\Psi^{(\bm 0, \bv_k)})} - \mathcal{T}^{(\Psi^{(\by+\bv_k,\by)})}}_F^2 \Bigg ),
\end{equation}
where $\norm{\cdot}_F$ denotes the Frobenius norm.

\subsubsection{Efficient evaluation of Frobenius norms for sums of low-rank tensors}

Let $k \in \{1,\dots,K\}$ and $\by \in S_M$ such that $\by + \bv_k \neq \bm 0$ be fixed. 
We derive an algorithm inspired by tensor train (TT) orthogonalization~\cite{Oseledets11} to evaluate $\norm{  \mathcal{T}^{(\Psi^{(\bm 0, \bv_k)})} - \mathcal{T}^{(\Psi^{(\by+\bv_k,\by)})}}_F^2$. 

We start by rewriting Equation~\eqref{eq:CPformat} as
\begin{equation}
    \label{eq:TTformat}
    T^{(\Psi)}_{i_1,\dots,i_N} = \sum_{k_1 = 1}^r \dots \sum_{k_{N-1} = 1}^r \mathcal{C}^{1}_{1,i_1,k_1} \mathcal{C}^{2}_{k_1,i_2,k_2} \mathcal{C}^{3}_{k_2,i_3,k_3} \cdots \mathcal{C}^{N}_{k_{N-1},i_N,1},   
\end{equation}
with tensors $\mathcal{C}^{1} \in \R^{1 \times  2 \times r}$, $\mathcal{C}^{N} \in \R^{r\times 2  \times 1}$ and $\mathcal{C}^{j} \in \R^{r\times 2 \times r}$ for $1 < j < N$, whose entries are given by
\begin{align*}
    \mathcal{C}^{j}_{k_1,i,k_2} &=  \left\{
\begin{array}{ll}
\bm a_i^{j,k_1} & k_1 = k_2 \\
0 & \, \textrm{otherwise} \\
\end{array}
\right., & \text{for } 1 < j < N,\ 1\leq i \leq 2,\ 1 \leq k_1,k_2 \leq N,  \\
\mathcal{C}^1_{1,i,k} &= \bm a_i^{1,k} & \text{for } 1\leq i \leq 2,\ 1 \leq k \leq N,   \\
\mathcal{C}^N_{k,i,1} &= \bm a_i^{N,k} & \text{for } 1\leq i \leq 2,\ 1 \leq k \leq N. 
\end{align*}

The representation format~\eqref{eq:TTformat} can be generalized to other low-rank tensors. Let $\mathcal{T}^{(\Psi^{(\by+\bv_k,\by)})}$ be analogously represented in the format~\eqref{eq:TTformat} by tensors $\mathcal{D}^j$. The tensor $\mathcal{T}^{(\Psi^{(\bm 0, \bv_k)})} - \mathcal{T}^{(\Psi^{(\by+\bv_k,\by)})}$ is at most rank-$2r$. It's representation in the form of~\eqref{eq:TTformat} with tensors $\mathcal E^j$ can be constructed from the tensors $\mathcal{C}^j,\mathcal D^j$, i.e. the tensors $\mathcal E^1 \in \R^{1 \times 2 \times 2r},\ \mathcal E^N \in \R^{2r \times 2 \times 1}$ and $\mathcal E^j \in \R^{2r \times 2 \times 2r}$ for $1<j<N$ are defined as $\mathcal{E}^1_{1,:,1:r} = \mathcal C^1,\ \mathcal{E}^1_{1,:,r+1:2r} = - \mathcal D^1,\  \mathcal{E}^N_{1:r,:,1} = \mathcal C^N,\ \mathcal{E}^N_{r+1:2r,:,1} = \mathcal D^N$ and $\mathcal{E}^j_{1:r,:,1:r} = \mathcal C^j,\ \mathcal{E}^j_{r+1:2r,:,r+1:2r} = \mathcal D^j$. Note that for $1<j<N$ these tensors are sparse by construction.

We can compute the norm $\norm{  \mathcal{T}^{(\Psi^{(\bm 0, \bv_k)})} - \mathcal{T}^{(\Psi^{(\by+\bv_k,\by)})}}_F^2$ directly form the tensors $\mathcal E^j$ in a fast and stable way using using TT orthogonalization~\cite{Oseledets11} as defined in Algorithm~\ref{alg:computeFrobNorm}. We want to emphasize that the sparsity structure of $\mathcal{E}^j$ is preserved in Algorithm~\ref{alg:computeFrobNorm}. The evaluation of $\norm{\mathcal{E}^{N}}_F^2$ in line~\ref{line:FrobNormEval} requires summing up $4r$ terms, whereas a naive evaluation of $\norm{  \mathcal{T}^{(\Psi^{(\bm 0, \bv_k)})} - \mathcal{T}^{(\Psi^{(\by+\bv_k,\by)})}}_F^2$ involves $2^N$ terms.

\begin{algorithm}
\caption{Frobenius norm evaluation}\label{alg:computeFrobNorm}
\begin{algorithmic}[1]
\algdef{SE}[SUBALG]{Indent}{EndIndent}{}{\algorithmicend\ }%
\algtext*{Indent}
\algtext*{EndIndent}
\State \textbf{Input}: tensors $\mathcal E^1 \in \R^{1 \times 2 \times 2r},\ \mathcal E^N \in \R^{2r \times 2 \times 1}$ and $\mathcal E^j \in \R^{2r \times 2 \times 2r}$ for $1<j<N$
\State \textbf{Output}: $\norm{\mathcal{T}}_F^2$, where $\mathcal T \in \R^{2^N}$ has entries 
$T_{i_1,\dots,i_N} = \sum_{k_1 = 1}^r \dots \sum_{k_{N-1} = 1}^r \mathcal{E}^{1}_{1,i_1,k_1} \mathcal{E}^{2}_{k_1,i_2,k_2} \cdots \mathcal{E}^{N}_{k_{N-1},i_N,1}$
\State Compute $QR$ decomposition of $\mathcal{E}^1$ reshaped as element in $\R^{2 \times 2r}$.
\State Set $\mathcal{E}^1$ to $Q$ reshaped as element in $\R^{1 \times 2 \times 2r}$.
\State Update $\mathcal E^2$: reshape to $\R^{2r \times 2\cdot 2r}$, multiply with $R$ from the left, reshape back to $\R^{2r \times 2 \times 2r}$.
\For $j = 2,\dots,N-1$
\State Compute $QR$ decomposition of $\mathcal{E}^j$ reshaped as element in $\R^{2r \cdot 2 \times 2r}$.
\State Set $\mathcal{E}^j$ to $Q$ reshaped as element in $\R^{2r \times 2 \times 2r}$.
\If $j < N-1$
\State Update $\mathcal E^{j+1}$: reshape to $\R^{2r \times 2\cdot 2r}$, multiply with $R$ from the left, reshape back to $\R^{2r \times 2 \times 2r}$.
\Else 
\State Update $\mathcal E^{N}$: reshape to $\R^{2r \times 2}$, multiply with $R$ from the left, reshape back to $\R^{2r \times 2 \times 1}$.
\EndIf
\EndFor
\State $\norm{  \mathcal{T}}_F^2 = \norm{\mathcal{E}^{N}}_F^2$ \label{line:FrobNormEval}
\end{algorithmic}
\end{algorithm} 

\begin{rmrk} 
The sum of several low-rank tensors can be represented in the format~\eqref{eq:TTformat}. The representation of the sum can again be constructed from the representations of the individual low-rank tensors. This allows us to apply TT orthogonalization to evaluate all Frobenius-norms in Equation~\eqref{eq:tensorformulation} in a fast and stable way.
\end{rmrk}

\subsection{Successive minimization}
Let $r \in \mathbb N$. In order to approximate the solution of the minimization problem  $\umin{\Phi \in \mathcal T_M^r} A_M^{\bm u}(\Phi)$, we decompose the problem into a sequence of successive minimization problems~\cite{Ammar12,Grasedyck13,Nouy17b}.
Let $\Phi \in H_M^r$ be represented as $\Phi = R^{(1)} + \dots + R^{(r)}$ by rank-$1$ functions $R^{(k)} \in H_M^1$. 
The main idea of successive minimization is to first determine $R^{(1)}$ as solution of $\min_{R^{(1)} \in \mathcal T_M^1} A_M^{\bm u}(R^{(1)})$. 
In a successive step $R^{(2)}$ is determined as solution of $\min_{R^{(2)} \in \mathcal T_M^1} A_M^{\bm u}(R^{(1)}+R^{(2)})$. 
This is continued successively until $\Phi$ is determined. 
This idea is formalized in Algorithm~\ref{alg:sucMin}.

\begin{algorithm}
\caption{Successive minimization}\label{alg:sucMin}
\begin{algorithmic}[1]
\algdef{SE}[SUBALG]{Indent}{EndIndent}{}{\algorithmicend\ }%
\algtext*{Indent}
\algtext*{EndIndent}
\State \textbf{Input}: rank $r$
\State \textbf{Output}: approximation $\Phi \approx \argmin_{\Phi \in \mathcal T_M^r} A_M^{\bm u}(\Phi)$ 
\State $\Phi = 0$
\For $k = 1,\dots,r$
\State $R^{(k)} =  \argmin_{R \in \mathcal T_M^1} A_M^{\bm u}(\Phi + R)$ 
\State $\Phi = \sum_{i = 1}^k  R^{(i)}$
\EndFor
\end{algorithmic}
\end{algorithm} 

\subsection{Alternating least squares}

In the successive minimization algorithm~\ref{alg:sucMin}, we need to solve minimization problems of the form  
\begin{equation}
    \label{eq:sucForm}
    \umin{R \in \mathcal T_M^1} A_M^{\bm u}(\Phi + R)
\end{equation}
for given $\Phi \in H_M^r$. 
In the following, we introduce an alternating least squares algorithm~\cite{Carroll70,Kolda09, de2008decompositions, oseledets2018alternating} to solve such minimization problems.
The main idea is to approximate the solution of~\eqref{eq:sucForm} by an iterative scheme which amounts to solving a sequence of small-dimensional linear problems. 
We start from an initial $R(\bEta) := \Pi_{\bs\in S_M}R_\bs(\eta_\bs)$. 
We first minimize $A^{\bm u}(\Phi + R)$ only with respect to a selected $R_{\bs_0}:\{0,1\} \to \R$ for some $\bs_0\in S_M$ leaving the other $R_\bs$, $\bs\neq \bs_0$ fixed. 
By partially evaluating  $A^{\bm u}(\Phi + R)$ for all terms not depending on $R_{\bs_0}$, we obtain that $\min_{R_{\bs_0} \in \{\{0,1\} \to \R\}}A^{\bm u} (\Phi + R)$ with $R(\bEta):= R_{\bs_0}(\eta_{\bs_0}) \Pi_{\bs\in S_M\setminus \{\bs_0\}} R_\bs(\eta_\bs)$ is equivalent to a quadratic optimization problem
\begin{equation} \label{eq:InterpolationPolynomial}
   \min_{R_{\bs_0} \in \{\{0,1\} \to \R\}} \alpha_1 R_{\bs_0}(1)^2 + \alpha_2 R_{\bs_0}(0)^2 + \alpha_3 R_{\bs_0}(1)R_{\bs_0}(0) + \alpha_4 R_{\bs_0}(1) + \alpha_5 R_{\bs_0}(0) + \alpha_6,
\end{equation}
with constants $\alpha_1,\dots,\alpha_6 \in \R$ depending on the fixed $R_\bs$, $\bs\neq \bs_0$ and $\Phi$. 
This quadratic optimization problem always admits a unique optimal $R_{s_0}$, which is given by $R_{\bs_0}(1) = a$ and $R_{\bs_0}(0) = b$, where $a,b\in R$ are the solution of the linear system 
\begin{equation}
    \begin{pmatrix}
2\alpha_1 & \alpha_3 \\ \alpha_3 & 2 \alpha_2
\end{pmatrix}
\begin{pmatrix}
a \\ b
\end{pmatrix} =
\begin{pmatrix}
-\alpha_4 \\ - \alpha_5
\end{pmatrix}.
\label{eq:Interpolation2x2System}
\end{equation}

This allows us to optimize $A^{\bm u}_M(\Phi + R)$ with respect to individual $R_{\bs_0}$. 
By alternating the the selected $\bs_0\in S_M$, we obtain the alternating least squares algorithm, which is formalized in Algorithm~\ref{alg:ALS}. 

\begin{algorithm}
\caption{Alternating least squares}\label{alg:ALS}
\begin{algorithmic}[1]
\algdef{SE}[SUBALG]{Indent}{EndIndent}{}{\algorithmicend\ }%
\algtext*{Indent}
\algtext*{EndIndent}
\State \textbf{Input}: initial functions $R^0_\bs: \{0,1\} \to \mathbb{R}$ for $\bs\in S_M$, function $\Phi \in H_M^r$, vector $\bm u$, tolerance $\eps$
\State \textbf{Output}: approximation ${R}_{\rm opt}(\bEta) = \Pi_{\bs\in S_M}R_\bs(\eta_\bs)$ of $\argmin_{R \in \mathcal T_M^1} A^u_M(\Phi + R)$
\State $v_{\textsf{old}} = \infty$, $v_{\textsf{new}} = A_M^{\bm u}(\Phi + R^0)$ with $R^0(\bEta):= \Pi_{\bs\in S_M} R^0_\bs(\eta_\bs)$
\State $\forall \bs\in S_M, \; R_\bs := R_\bs^0$. 
\While $\abs{v_{\textsf{old}}-v_{\textsf{new}}} > \eps \abs{v_{\textsf{new}}}$.
\State $v_{\textsf{old}} = v_{\textsf{new}}$
\For $\bs_0\in S_M$
\State $R_{\bs_0} = \argmin_{\widetilde{R}_{\bs_0} :\{0,1\} \to \R} A^{\bm u}(\Phi + \widetilde{R})$ 
where $\widetilde{R}(\bEta) = \widetilde{R}_{\bs_0}(\eta_{\bs_0}) \Pi_{\bs\in S_M \setminus \{\bs_0\}} R_\bs(\eta_\bs)$ for all $\bEta = (\eta_\bs)_{\bs\in S_M}$ \label{line:optimizesinglefunction}
\EndFor
\State $v_{\textsf{new}} = A_M^{\bm u}(\Phi + R)$
\EndWhile
\end{algorithmic}
\end{algorithm} 

\begin{rmrk}
To compute the constants $\alpha_i$, we can either explicitly implement the partial evaluations of $A^{\bm u}_M(\Phi + R)$. 
Alternatively, we can treat $A^{\bm u}_M(\Phi + R)$ as a function in $\R^2\to \R$ depending on the values $R_{\bs_0}(0)$ and $R_{\bs_0}(1)$. We know that this function is a multivariate-polynomial of the form $\alpha_1 R_{\bs_0}(1)^2 + \alpha_2 R_{\bs_0}(0)^2 + \alpha_3 R_{\bs_0}(1)R_j(0) + \alpha_4 R_{\bs_0}(1) + \alpha_5 R_{\bs_0}(0) + \alpha_6$. The constants can be computed using  multivariate-polynomial interpolation in six points. This interpolation has the advantages that it is non-intrusive and that the evaluations of $A^{\bm u}_M(\Phi+R)$ can be performed efficiently using the ideas of Section~\ref{sec:StableEvaluation}.
\end{rmrk}

\begin{rmrk}
Throughout this work, we use approximations in the CP-format. The matrix product state representation~\cite{Perez07} also known as tensor train format~\cite{Oseledets11} would present an alternative low-rank format, for which minimization problems can also be solved using alternating algorithms~\cite{Holtz12,Grasedyck13,Szalay15}. 
\end{rmrk}

\subsection{Monte Carlo methods}\label{sec:MonteCarlo}

Let $\Phi \in H^r_M$ denote an approximation of the solution of~\eqref{eq:minopt}. 
In the following, we discuss how to evaluate $A_{M,\ell}^{\bm u}(\Phi)$. 
In the definition of $A_{M,\ell}^{\bm u}$ the function 
\begin{equation}
    \label{eq:RelevantSamplingTerm}
    f(\bEta) := \sum_{k=1}^K p_k
\Bigg( (1- \eta_{\bv_k}) \bra{
\bu\cdot \bv_k + \Psi(\bEta^{0,\bv_k}) - \Psi(\bEta)}^2 + \frac{1}{2}
\sum_{\substack{\by \in S \setminus \{ {\bm 0} \} \\ \by+\bv_k \neq \bm 0 } } \bra{
\Psi(\bEta^{\by+\bv_k,\by})-\Psi(\bEta)}^2   \Bigg)
\end{equation}
is evaluated for all $\bEta \in C_{M,\ell}$.
For larger $N$ and most choices of $\ell$ evaluating $|C_{M,\ell}| = \binom{N}{\ell}$ terms is intractable. 
We thus propose to use a Monte Carlo method~\cite{Metropolis49} to approximate $A_{M,\ell}^{\bm u}(\Phi)$.

In a naive Monte Carlo method we compute $N_s$ samples $\bEta^{(i)}\in C_{M,\ell}$ for $1\leq i \leq N_s$ and replace 
$\frac{1}{|C_{M,\ell}|} \sum_{\bEta \in C_{M,\ell}}$ by $\frac{1}{N_s} \sum_{\bEta^{(1)},\dots,\bEta^{({N_s})}}$. In Algorithm~\ref{alg:varianceReduc}, we define a Monte Carlo method with additional variance reduction~\cite{Hammersley65}, which as demonstrated in Section~\ref{sec:numerical:experiments} reduces the number of samples needed to obtain a given approximation accuracy. The main idea is to observe that the sites $\bm v_1,\dots, \bm v_K$ play a special role since they are the most relevant for jumps of the tagged particle. Instead of sampling the $\bEta$ uniformly in $C_{M,\ell}$, we now ensure that all possible states of the sites $\bm v_1,\dots, \bm v_K$ are occurring equally often in the set of sample points.

\begin{algorithm}
\caption{Monte Carlo method with variance reduction}\label{alg:varianceReduc}
\begin{algorithmic}[1]
\algdef{SE}[SUBALG]{Indent}{EndIndent}{}{\algorithmicend\ }%
\algtext*{Indent}
\algtext*{EndIndent}
\State \textbf{Input}: function $f:C_{M,\ell} \to \mathbb R$, parameter $\widetilde{N}_s$
\State \textbf{Output}: approximation $S$ of $\frac{1}{|C_{M,\ell}|}\sum_{\bEta \in C_{M,\ell}} f(\bEta)$
\State $S = 0$
\For $\vec{\bm w}^{(1)} \in \{0,1\}^K$ 
    \State $n_1 := || \vec{\bm w} ||_1$, $n_2 := \ell - n_1$
    \If $0 \leq n_2 \leq   N-\ell$
        \State $\Tilde{S}=0$
        \For $i = 1,\dots,\widetilde{N}_s$
            \State Sample $\vec{\bm w}^{(2)} \in \{0,1\}^{S_M\setminus \{ \bm v_1,\dots,\bm v_k\}}$ such that $|| \vec{\bm w}^{(2)}||_1 = n_2$. \label{line:sampleW2}
            \State Construct $\bEta \in C_{M,\ell}$ such that $\bEta(\bm s) = \left\{
\begin{array}{ll}
\vec{\bm w}^{(1)}_k & \bm s = \bm v_k \\
\vec{\bm w}^{(2)}_{\bm s} & \text{otherwise} \\
\end{array}
\right..$ 
        \State $\Tilde{S} = \Tilde{S} + f(\bEta)$
        \EndFor        
        \State $S = S + \binom{N-K}{n_2}  \cdot \Tilde{S} /\widetilde{N}_s$ \label{line:SumS1}
    \EndIf
\EndFor
\State $S = \frac{S}{\binom{N}{\ell}}$ \label{line:SumS2}
\end{algorithmic}
\end{algorithm} 

\begin{rmrk}
In line~\ref{line:sampleW2} of Algorithm~\ref{alg:varianceReduc}, we sample form the set set $\{\vec{\bm w}^{(2)} \in \{0,1\}^{S_M\setminus \{ \bm v_1,\dots,\bm v_k\}}  | || \vec{\bm w}^{(2)}||_1 = n_2 \}$ which contains $\binom{N-K}{n_2}$ elements. When the number of elements in this set is smaller than $\widetilde{N}_s$, we can compute $\Tilde{S}$ based on all possible $\vec{\bm w}^{(2)}$ instead of sampling $\widetilde{N}_s$ times. This decreases the variance of the approximation further and reduces the number of required evaluations of $f$. 
\end{rmrk}

\begin{rmrk}
Lines~\ref{line:SumS1} and~\ref{line:SumS2} of Algorithm~\ref{alg:varianceReduc} potentially lead to stability issues in floating point arithmetic when $\binom{N-K}{n_2}$ or $\binom{N}{\ell}$ are large. This problem can be mitigated using the log-sum-exp trick.
\end{rmrk}

\subsection{\corr{Limitations of the approach}}
\corr{In this section, we briefly discuss the limitations of the proposed algorithm. These are primarily related to the following observation.
Since $A_{M,\ell}^{\bm u}(\Psi^{M, \bm u}_{\textsf{opt}})$ defines the entries of the self-diffusion coefficient~\eqref{eq:approx:Ds}, we know that 
$A_{M,\ell}^{ \bm u}(\Psi^{M, \bm u}_{\textsf{opt}}) \in [0,1]$.
Note that $A_{M,\ell}^{\bm u}$ contains the normalization constant $|C_{M,\ell}|^{-1}$. 
The objective function $A_{M}^{ u}$~\eqref{eq:def:Au} does not include a normalization factor.
This implies that the $A_{M}^{\bm u}(\Psi^{M, \bm u}_{\textsf{opt}})$ is of order $2^N$. 
Trying to minimize this objective function numerically using floating point arithmetic leads to issues caused by rounding errors for large $N$.}

\corr{For larger values of the objective function, it  becomes increasingly challenging to solve the individual minimization problems in line~\ref{line:optimizesinglefunction} of Algorithm~\ref{alg:ALS}.
In particular, our approach of using polynomial interpolation to find the location of the minimum might encounter numerical rounding issues.
We find that the positive eigenvalues of the $2\times 2$ matrix in~\eqref{eq:Interpolation2x2System} tend to be many orders of magnitude smaller than the norm of the right hand side. 
This implies that small rounding errors in the constants $\alpha_1,\dots,\alpha_6$ might lead to vastly different solutions.
For $M>2$, we even find that rounding errors can lead to negative eigenvalues in the system, i.e. we can not find the minimum of the  polynomial~\eqref{eq:InterpolationPolynomial} at all.
It might be necessary to use a different approach to solve the individual ALS minimization problems for larger $M$.
Alternatively, one could develop an algorithm to minimize the function $\log(A_{M}^{\bm u})$.
}

\corr{Moreover, the number of ALS iterations needed to find a good rank-$1$ minimizer tends to increase when increasing the size of the domain. 
This is especially true when a random initialization is used in Algorithm~\ref{alg:ALS}.
Note that it might be possible to circumvent this issue by initializing Algorithm~\ref{alg:ALS} based on the solution computed for a smaller value of $M$.
Overcoming these limitations will then require further investigation which we intend to carry out in a following work. A possible path could be to combine domain decomposition approaches together with the tensor-based optimization algorithm we propose here, in order to obtain a global optimization procedure where only independent parallel local optimization problems on medium-sized cells are solved.
}

\section{Numerical Experiments}\label{sec:numerical:experiments}

In the following numerical experiments~\footnote{{The code to reproduce all experiments is available on \url{https://github.com/cstroessner/SelfDiffusion}}}, we consider the a jumping model with $K=4$ displacement vectors $\bv_1=(1,0)$, $\bv_2=(-1,0)$, $\bv_3=(0, 1)$, and $\bv_4=(0,-1)$ and with associated probability $p_k=1/4$ for $k=1,2,3,4$.
All computing times are measured without parallelization in MATLAB R2018b on a Lenovo Thinkpad T480s with Intel Core i7-8650U CPU and 15.4 GiB RAM. Note that both the sampling Algorithm~\ref{alg:StochasticAlgorithm} and Algorithm~\ref{alg:varianceReduc} as well as the evaluation of the different Frobenius-norms in Equation~\eqref{eq:tensorformulation} can be parallelized to speed up computations.

\subsection{Solving the optimization problem}
In the following, we analyze the numerical approximation of the self-diffusion coefficient by solving the optimization problem~\eqref{eq:approx:Ds}. 

\medskip

\paragraph{Low-rank approximation error}
First, we study how the rank $r$ affects $\min_{\Phi \in H_M^r} A_M^{\bm u}(\Phi)$ for $\bm u = (1,0)$. Let $\Phi_{\textsf{ALS}}^{r,\bm u} \in H_M^r$ denote the function obtained by using Algorithm~\ref{alg:sucMin} with rank $r$, where the minimization problem in each iteration is solved using Algorithm~\ref{alg:ALS}.  The initial functions $R_{\bm s}^0$ in Algorithm~\ref{alg:ALS} are selected randomly by assigning random values drawn form the uniform distribution on $[0,1]$ to $R_{\bm s}^0(0),R_{\bm s}^0(1)$. In Algorithm~\ref{alg:ALS}, we set $\eps = 10^{-12}$ and additionally stop after line~\ref{line:optimizesinglefunction} has been executed $420$ times.

\begin{figure}[!ht]
    \centering 
    \includegraphics[width = 0.45\textwidth]{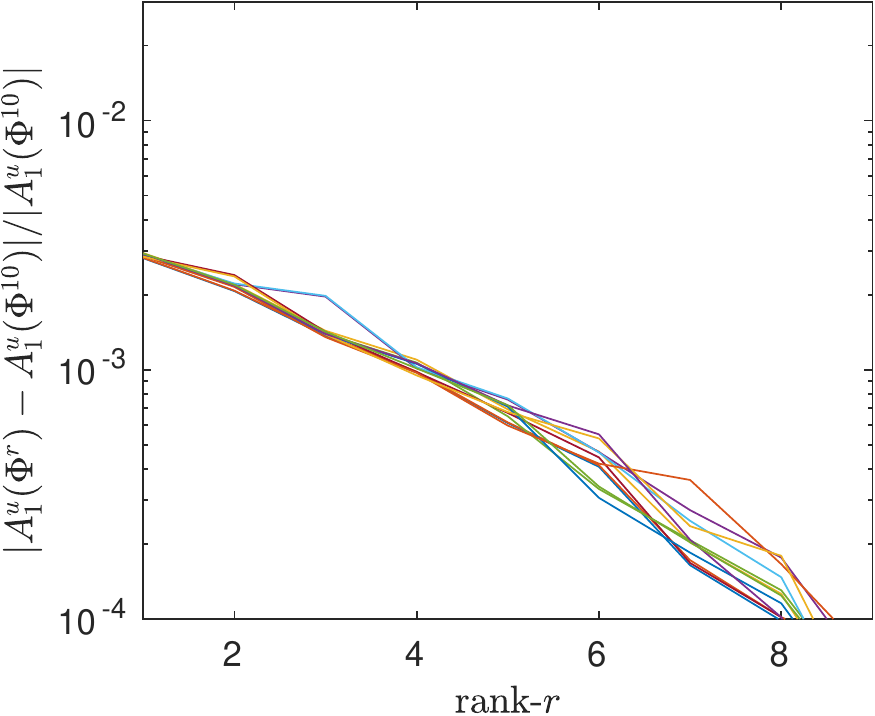} \hfill
    \includegraphics[width = 0.45\textwidth]{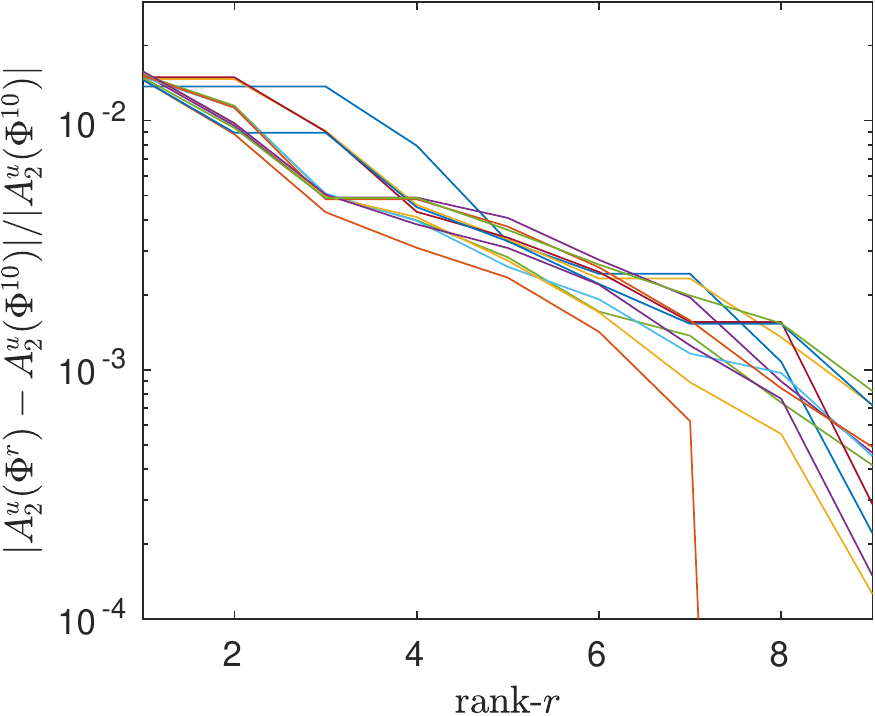}
    \caption{
    Algorithm~\ref{alg:sucMin} yields successive approximations $\Phi_{\textsf{ALS}}^{r,\bm u} \in H_M^r$. 
    For $r \in \{1,\dots,9\}$, we plot the relative error of $\Phi_{\textsf{ALS}}^{k,\bm u}$ compared to $\Phi_{\textsf{ALS}}^{10,\bm u}$.
    We repeat this experiment $12$ times with different random initial functions in Algorithm~\ref{alg:ALS}.
    The resulting relative errors are displayed in different colors.
    Left: $M=1$. Right: $M=2$.}
    \label{fig:RankImpact}
\end{figure}

We depict the error of $\Phi_{\textsf{ALS}}^{r,\bm u}$ for various $r$ in Figure~\ref{fig:RankImpact}. The computation with $r=10$ takes $4$ minutes for $M=1$ and $82$ minutes for $M=2$. We observe that the error decays quickly with increasing rank. In particular, for $M=1$ an approximation with $r=1$ is less than $0.1\%$ away from a direct least squares solution of the minimization problem~\eqref{eq:minopt}. Such a direct solution is intractable for $M=2$. 

\medskip

\begin{figure}[!ht] 
    \centering 
    \includegraphics[width =0.45\textwidth]{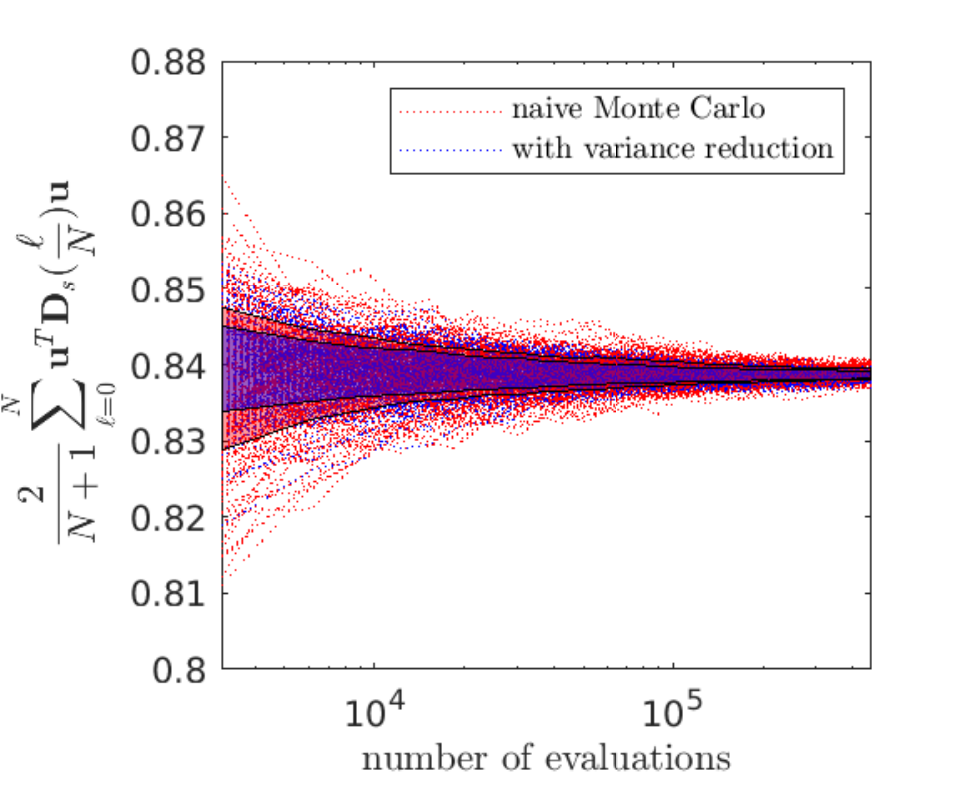} \hfill
    \includegraphics[width =0.45\textwidth]{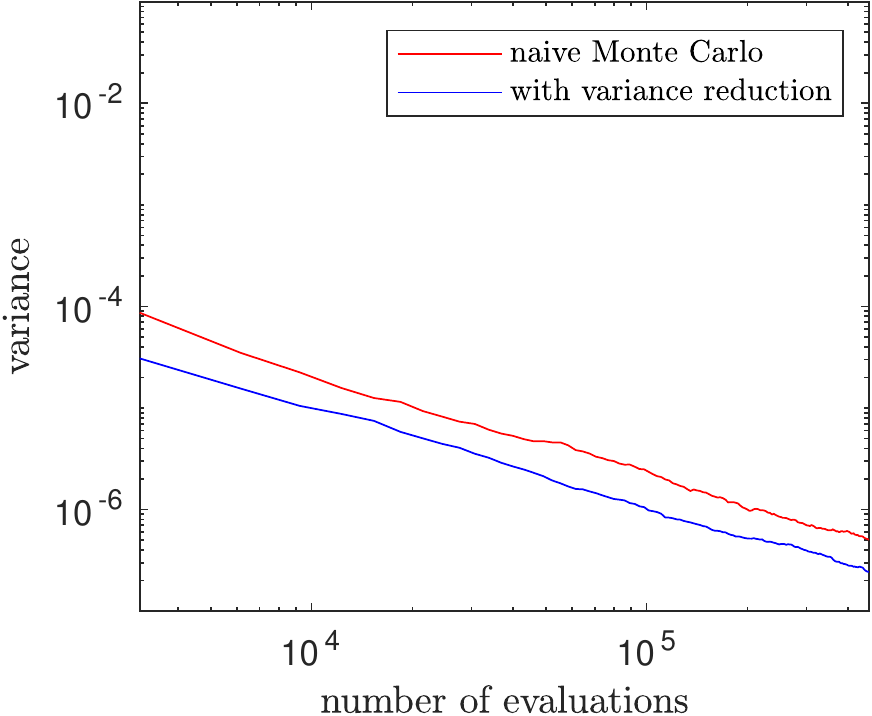}
    \caption{
    For $M=2$ and $\bm u=(1,0)$ we compute $\Phi_{\textsf{ALS}}^{\corr{3},\bm u}$ using Algorithm~\ref{alg:sucMin}. 
    We then approximate $\bm u^T \matD_s(\frac{\ell}{N}) \bm u = 2 A^{\bm u}_{2,\ell}(\Phi_{\textsf{ALS}}^{\corr{3},\bm u})$ for $0\leq \ell \leq N$ using sampling as in Section~\ref{sec:MonteCarlo}.
    We sample $\bEta \in C_{M,\ell}$ using two different approaches; from the uniform distribution (naive Monte Carlo) and using Algorithm~\ref{alg:varianceReduc} (with variance reduction).
    We use a total of $460\,800$ evaluations of~\eqref{eq:RelevantSamplingTerm} for both sampling methods. 
    After every $3\,072$ evaluations, we evaluate $\frac{2}{N+1} \sum_{\ell = 0}^N \bm u^T \matD_s(\frac{\ell}{N}) \bm u$  based on the current approximation.
    The whole experiment is repeated $250$ times.
    Left: 
    Evolution of $\frac{2}{N+1} \sum_{\ell = 0}^N \bm u^T \matD_s(\frac{\ell}{N}) \bm u$ with increasing number of samples. 
    The shaded areas mark one standard deviation from the mean for the respective algorithm. 
    Right: 
    Evolution of the variance of $\frac{2}{N+1} \sum_{\ell = 0}^N \bm u^T \matD_s(\frac{\ell}{N}) \bm u$ with increasing number of samples.}
    \label{fig:OptimVariance}
\end{figure}

\paragraph{Sampling-based evaluation} 
Let $M=2$ and $\bm u = (1,0)$. For a given approximate solution $\Phi_{\textsf{ALS}}^{\corr{3},\bm u}$, we want to evaluate $A^{\bm u}_{2,\ell}(\Phi_{\textsf{ALS}}^{\corr{3},\bm u})$.
A single evaluation of~\eqref{eq:RelevantSamplingTerm} with $\Psi$ replaced by $\Phi_{\textsf{ALS}}^{\corr{3},\bm u}$ requires on average around $\corr{5.1 \cdot 10^{-4}}$ seconds. 
Computing $A^{\bm u}_{2,\ell}(\Phi_{\textsf{ALS}}^{\corr{3},\bm u})$ for $0\leq \ell \leq N$ directly would require $2^{24}\approx 1.6\cdot 10^{7}$ evaluations of~\eqref{eq:RelevantSamplingTerm}, i.e. around \corr{$2.4$} hours. 
In Section~\ref{sec:MonteCarlo}, we proposed two Monte Carlo algorithms to obtain approximations of $A^{\bm u}_{2,\ell}(\Phi_{\textsf{ALS}}^{\corr{3},\bm u})$.
We visualize the variance in the approximation obtained by these algorithms in Figure~\ref{fig:OptimVariance}. 
The studied quantity $\frac{2}{N+1} \sum_{\ell = 0}^N \bm u^T \matD_s(\frac{\ell}{N}) \bm u$ approximates $\int_0^1 {\rm Tr} (\matD_s(\overline{\rho})) d\overline{\rho}$, where ${\rm Tr}$ denotes the trace operator.
Algorithm~\ref{alg:varianceReduc} clearly leads to a variance reduction.
With only $10^5$ samples, which can be evaluated in about one minute, we can already reach a variance of $10^{-6}$.

\subsection{Estimation of long-time mean square deviation}
In the following, we study estimating the long-term limit using Algorithm~\ref{alg:StochasticAlgorithm}.

\begin{figure}[!ht] 
    \centering 
    \includegraphics[width = 0.45\textwidth]{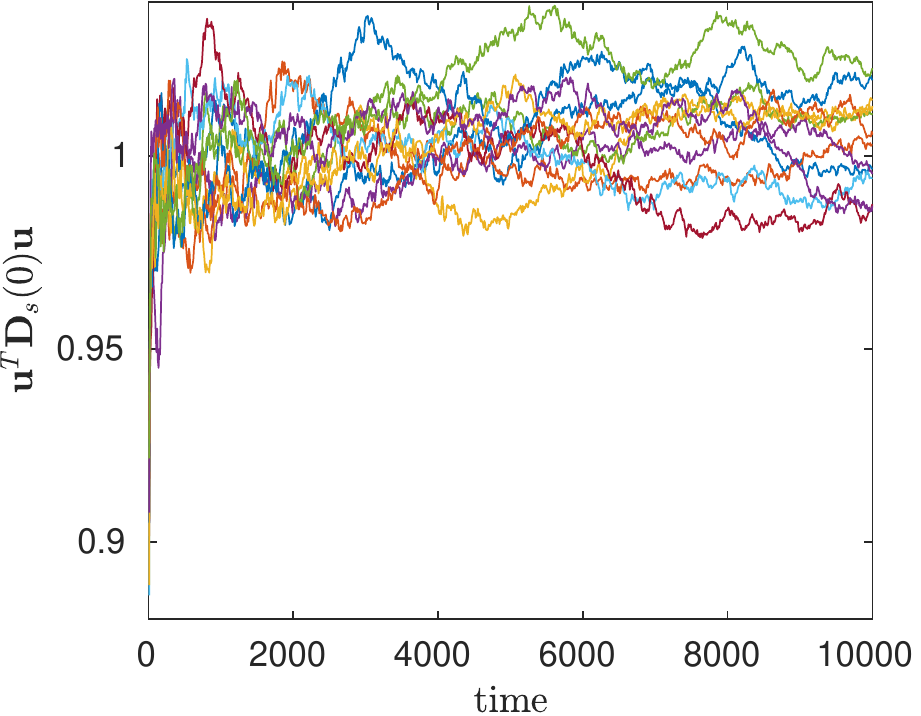} \hfill
    \includegraphics[width = 0.45\textwidth]{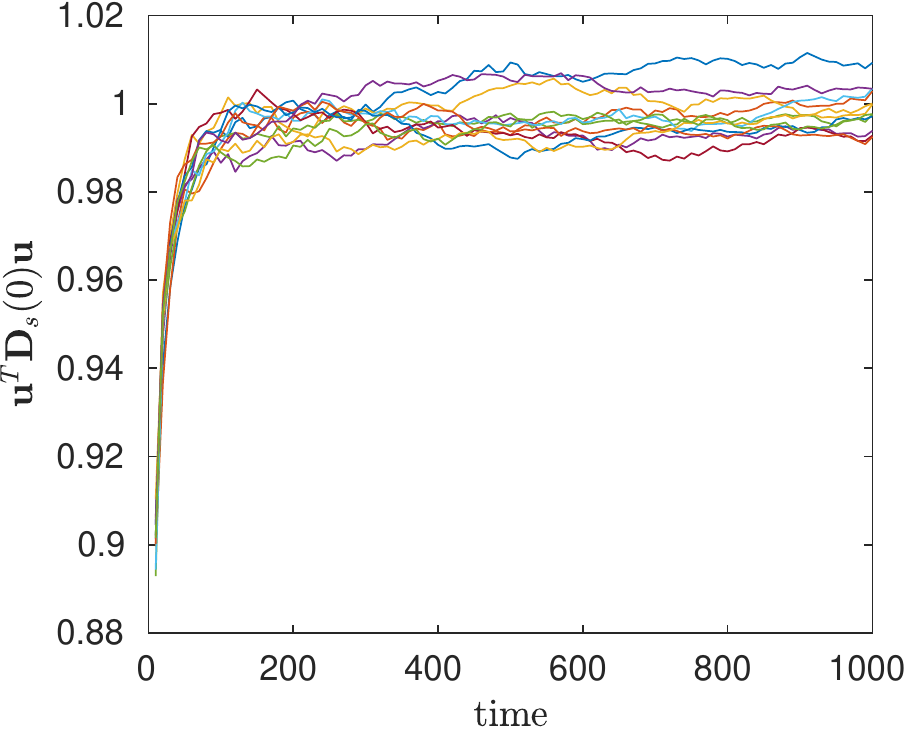}
    \caption{We run Algorithm~\ref{alg:StochasticAlgorithm} with $M=2$, $\bm u = (1,0)$, $\ell = 0$ and different values for $T$ and $\hat N_s$. For each setting, we plot how the approximation of the value $\bm u^T \matD_s(0) \bm u$, which is known to be equal to $1$, evolves over time for $12$ random initializations visualized in different colors. Left: $\hat N_s = 10\,000$, $T=10\,0000$. Right: $\hat N_s = 100\,000$, $T=1\,000$.}
    \label{fig:CalibratingTandNs}
\end{figure}

Figure~\ref{fig:CalibratingTandNs} motivates our choice of  $T=300$ and $\hat{N}_s = 30\,000$ for the following numerical experiments. The error caused by stopping with $T=300$ is negligible compared to the stochastic variance when stopping with  $\hat N_s = 30\,000$.

\begin{figure}[!ht] 
    \centering 
    \includegraphics[width = 0.45\textwidth]{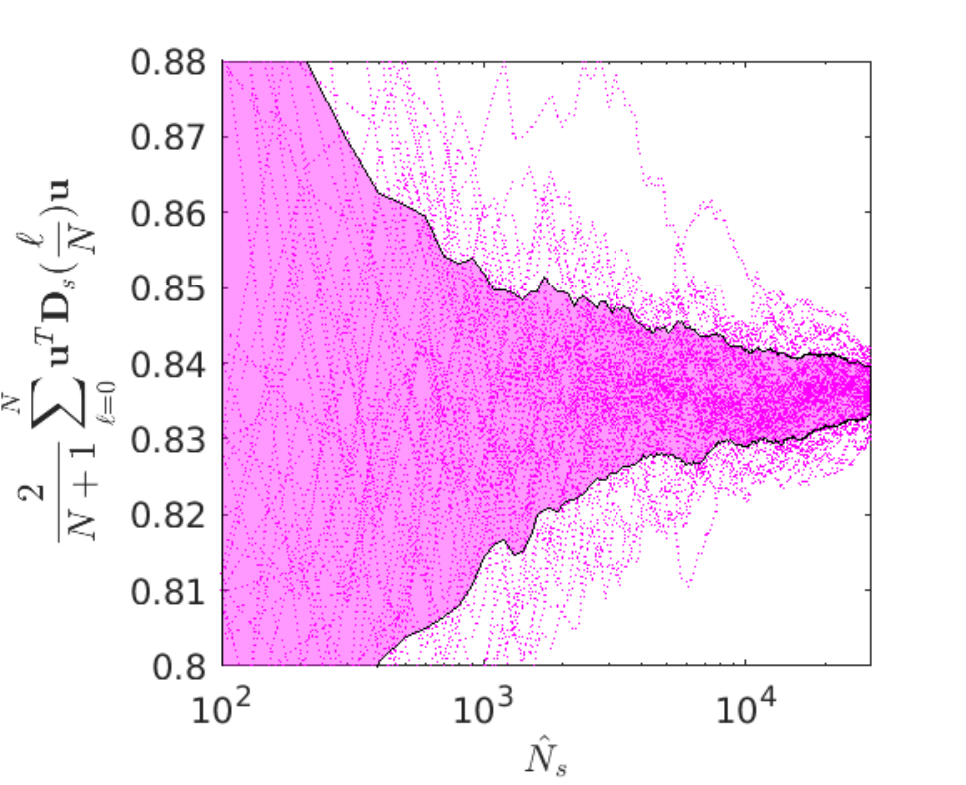}\hfill
    \includegraphics[width = 0.45\textwidth]{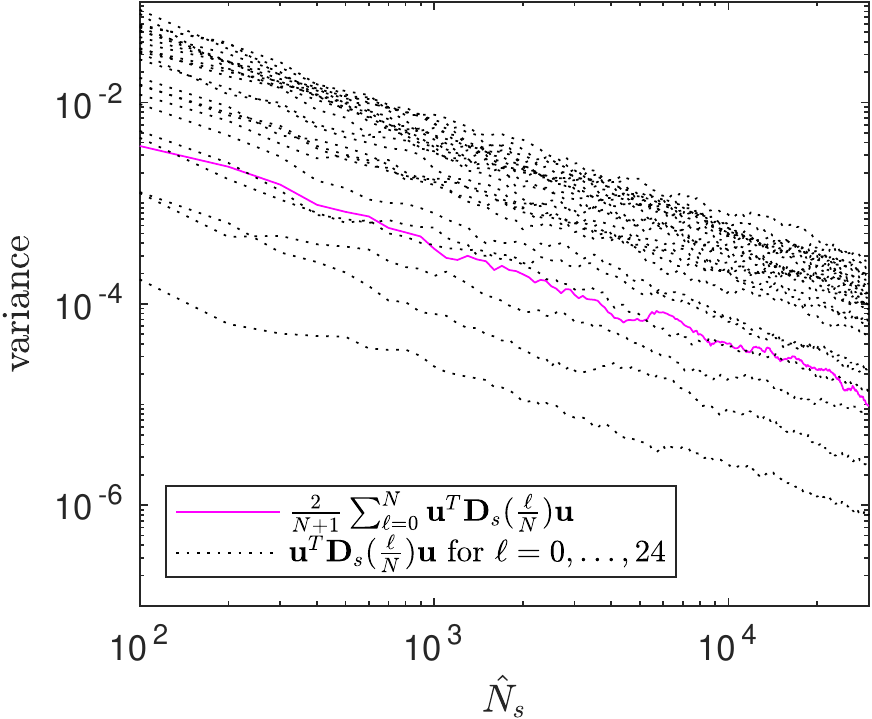}
    \caption{We run Algorithm~\ref{alg:StochasticAlgorithm} with $T=300$, $\bm u = (1,0)$ for $\ell = 0,\dots,24$ to approximate $\matD_s(\frac{\ell}{N})$. We use $\hat{N}_s = 30\,000$ and store intermediate values for $\hat{N}_s = 100,200,300,\dots$.
    Based on these approximations we compute $\frac{2}{N+1} \sum_{\ell = 0}^N \bm u^T \matD_s(\frac{\ell}{N}) \bm u$.
    This procedure is repeated $50$ times with different random initializations. 
    Left: Evolution of $\frac{2}{N+1} \sum_{\ell = 0}^N \bm u^T \matD_s(\frac{\ell}{N}) \bm u$ with increasing values of $\hat{N}_s$. The shaded area marks one standard deviation from the mean. 
    Right:  Evolution of the variance of $\frac{2}{N+1} \sum_{\ell = 0}^N \bm u^T \matD_s(\frac{\ell}{N}) \bm u$ with increasing values of $\hat{N}_s$. We additionally display the evolution of the variance of $ \bm u^T \matD_s(\frac{\ell}{N}) \bm u$ for individual values of $\ell$.}
    \label{fig:SamplingStochasicVariance}
\end{figure}

\begin{figure}[!ht] 
    \centering 
    \includegraphics[width = 0.45\textwidth]{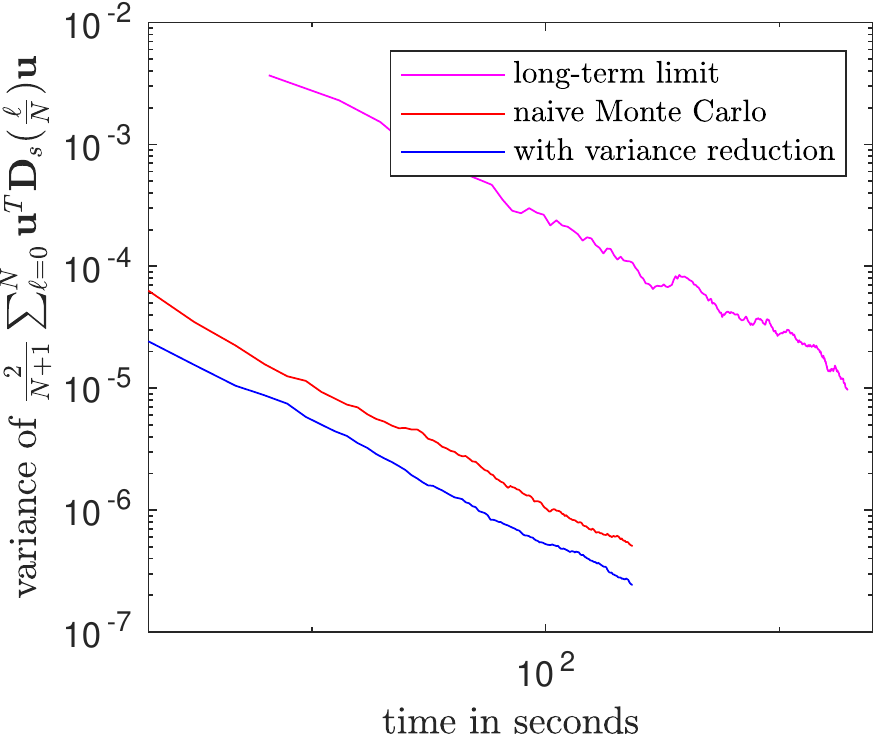}
    \hfill
    \includegraphics[width = 0.45\textwidth]{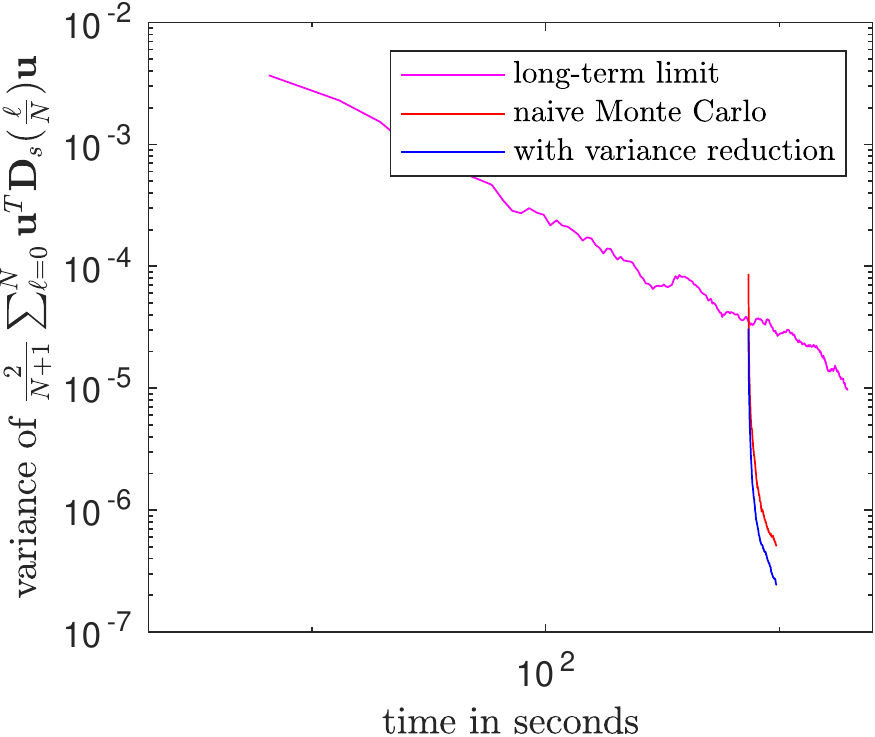}
    \caption{\corr{We plot the variances displayed on the right of Figure~\ref{fig:OptimVariance} and Figure~\ref{fig:SamplingStochasicVariance}.
    Instead of the number of samples, we display the computational time for evaluating the corresponding number of samples on the x-axis.
    Left: Computational times only for Algorithm~\ref{alg:StochasticAlgorithm}, the naive Monte Carlo method and Algorithm~\ref{alg:varianceReduc}.
    Right: We additionally include the computational time of $736$ seconds for computing $\Phi_{\textsf{ALS}}^{{3},\bm u}$ using Algorithm~\ref{alg:sucMin}.}
    }
    \label{fig:VarianceRuntime}
\end{figure}

In Figure~\ref{fig:SamplingStochasicVariance} we analyze how the parameter $\hat{N}_s$ affects an approximation of $\frac{2}{N+1} \sum_{\ell = 0}^N \bm u^T \matD_s(\frac{\ell}{N}) \bm u$. We observe that the quantity of interest converges to the same value around 0.84 as in Figure~\ref{fig:OptimVariance}. Note that a single run with $\hat{N}_s = 30\,000$ requires $33$ minutes and reaches a variance of $10^{-4}$. Compared to the results in Figure~\ref{fig:OptimVariance} the computational time needed to to achieve the same variance with Algorithm~\ref{alg:StochasticAlgorithm} compared to Algorithm~\ref{alg:varianceReduc} is much higher.
\corr{This is visualized in Figure~\ref{fig:VarianceRuntime}.}

\subsection{Comparison of algorithms}

We now compare the two approaches of approximating the self-diffusion coefficient. We evaluate $\bm u ^T \matD_s(\ell/N) \bm u$ for $\bm u \in \{ (1,0),(0,1),(1,1)\}$ to recover all entries of the symmetric matrix $\matD_s(\ell/N)$. Interpolation of the entries yields an approximation of $\matD_s(\rho)$ for $\rho \in [0,1]$. Alternatively, interpolation in the space of symmetric positive definite matrices could be used~\cite{Moakher06}.

For the following experiments, we use Algorithm~\ref{alg:StochasticAlgorithm} with a slight modification to approximate $\bu^T \matD_s(\overline{\rho}) \bu$ for different values of $\bm u$ simultaneously. This is achieved by keeping track of different $\alpha$ for different $\bm u$ in lines~\ref{line:UpdateAlphaBasedonU} and~\ref{line:UpdateAlpha2}. 
The optimization problem~\eqref{eq:minopt} needs to be solved for every $\bm u$ separately, but the solution can be evaluated for different $\ell$. 
In contrast, the sampling in Algorithm~\ref{alg:StochasticAlgorithm} needs to be redone completely for every $\ell$.

\begin{figure}[!ht]
    \centering
    \includegraphics[width = 0.45\textwidth]{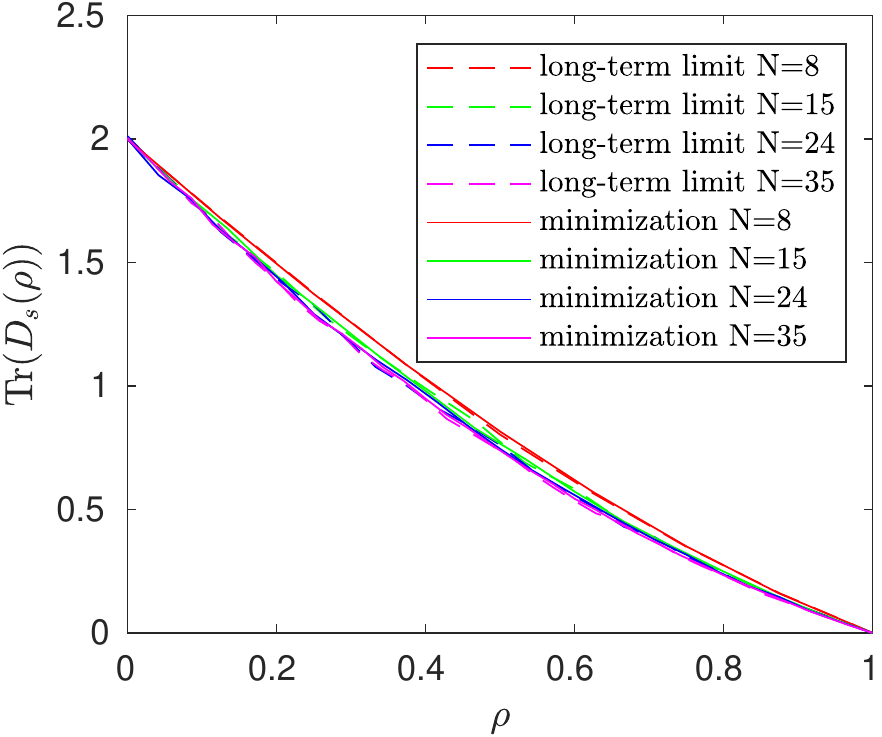}
    \hfill
    \includegraphics[width = 0.45\textwidth]{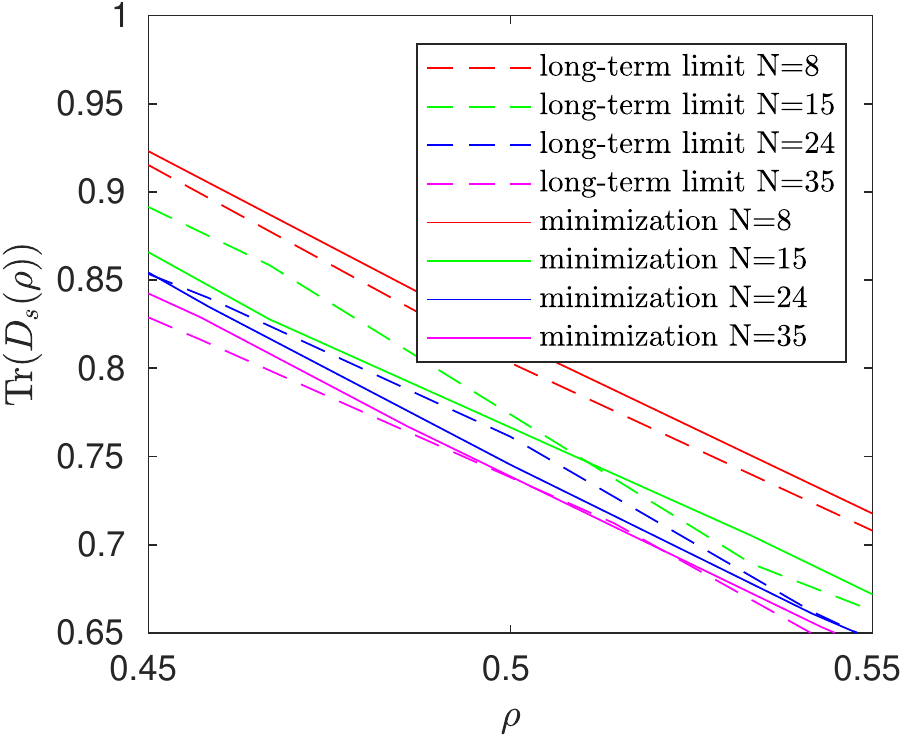}
    \caption{
    For $M=1$ ($N=8$) and $M=2$ ($N=24$) we solve the minimization problem~\eqref{eq:minopt} using Algorithm~\ref{alg:sucMin} with $r=3$.
    We evaluate $A^{\bm u}_{M,\ell}$ directly for $M=1$, whereas for $M=2$ we use Algorithm~\ref{alg:varianceReduc} with $\widetilde{N}_s = 50$.
    We compare to Algorithm~\ref{alg:StochasticAlgorithm} to estimate the long-term limit~\eqref{eq:infiniteLimit} with $\hat{N}_s = 30\,000$ and $T=300$ for $M=1$ and $M=2$ to compute $\matD_s(\ell/N)$.
    Repeating this for $\bm u \in \{ (1,0),(1,1),(0,1)\}$ yields $\matD_s(\ell/N)$.
    Element-wise linear interpolation allows us to evaluate $\matD_s(\rho)$ for $\rho \in [0,1]$.
    We plot the trace of $\matD_s(\rho)$.
    In addition, we run both algorithms on a $4\times 4$ \corr{and $6\times 6$} periodic grid ($N=15$ \corr{ and $N=35$}). \corr{For $N=15$ we use the same parameters as for $N=24$. For $N=35$ we use rank $10$, $\widetilde{N}_s = 150$ and $\hat{N}_s = 80\,000$, $T=600$.} Note that \corr{the $4\times 4$ grid lies} in between the $3 \times 3$ grid for $M=1$ and the $5 \times 5$ grid for $M=2$.
    Right: Zoom on part of the graphs.}
    \label{fig:comparisonToStochasticApproach}
\end{figure}

\begin{figure}[!ht]
    \centering
    \includegraphics[width = 0.32\textwidth]{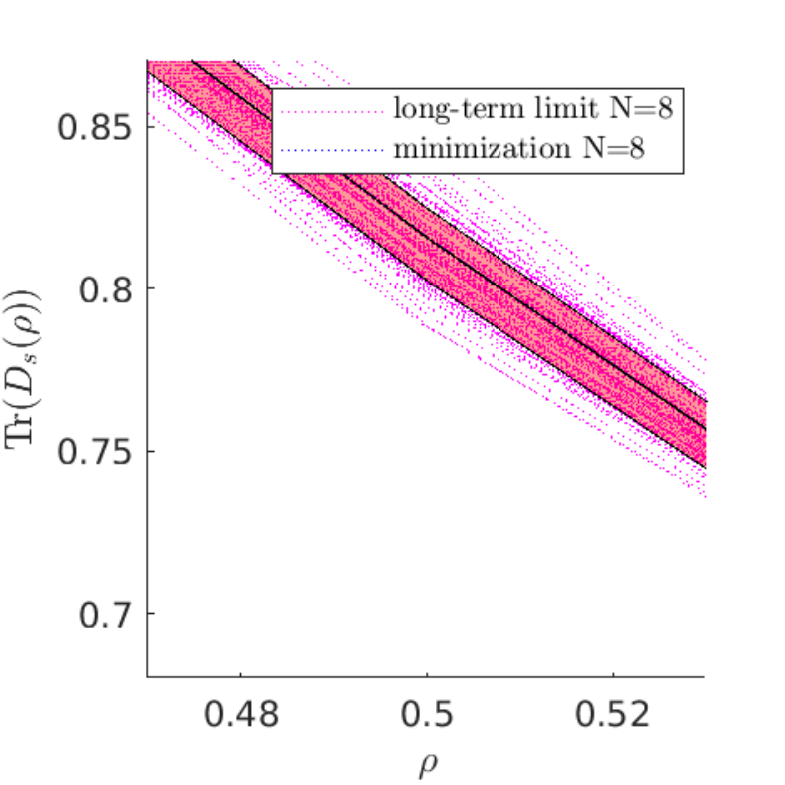}
    \hfill
    \includegraphics[width = 0.32\textwidth]{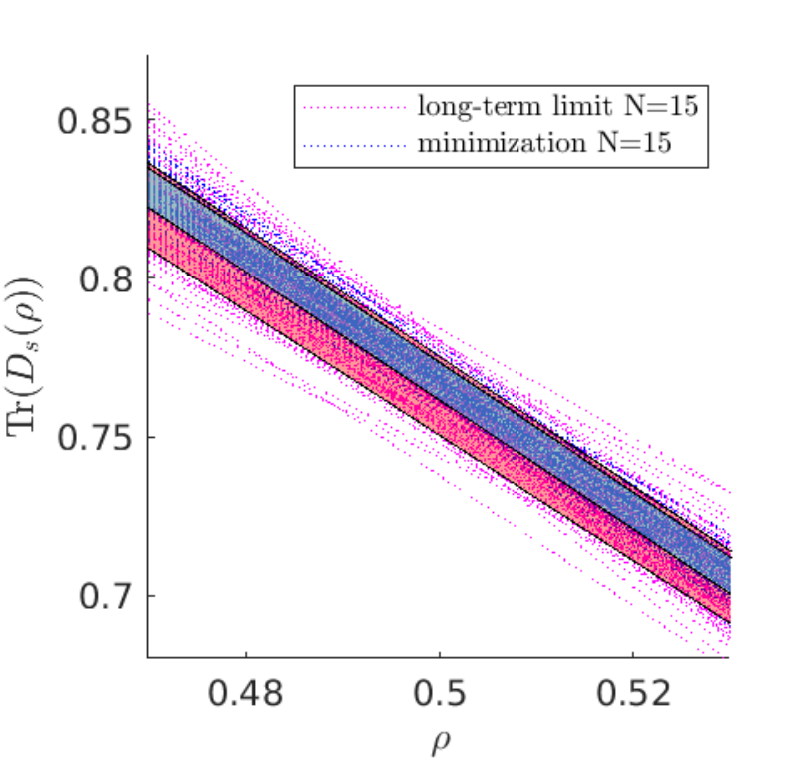}
    \hfill
    \includegraphics[width = 0.32\textwidth]{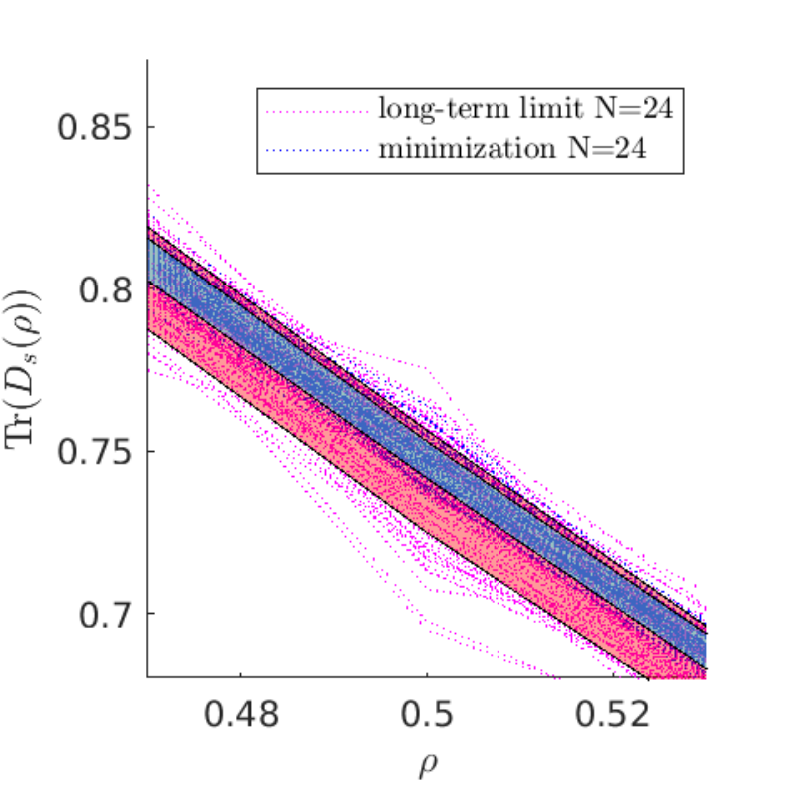}
    \caption{
    \corr{We repeat the experiment presented in Figure~\ref{fig:comparisonToStochasticApproach} $100$ times for $N=8$, $N=15$ and $N=24$. We plot a zoom on a part of the graphs. The shaded areas mark one standard deviation from the mean for each of the methods.}}
    \label{fig:traceVariance}
\end{figure}

In Figure~\ref{fig:comparisonToStochasticApproach}, we display the trace of $\matD_s(\rho)$ computed using different approaches.
Computing the whole matrix $\matD_s(\ell/N)$ for all $0\leq \ell \leq N$ took $3$ minutes ($37$ minutes) for the approach based on solving the minimization problem, whereas using Algorithm~\ref{alg:StochasticAlgorithm} took $15$ minutes ($45$ minutes) for $M=1$ ($M=2$). 
The much higher variance of Algorithm~\ref{alg:StochasticAlgorithm} for approximating the long-term limit leads to clearly visible changes in the graph of trace at the values $\rho = \frac{\ell}{N}$. \corr{The figure also contains results using a $6\times 6$ grid ($N=35$) for which the minimization based approach took 5 hours and 40 and the sampling of the long-term limit took 9 hours and 24 minutes.}
The sampled solution of the minimization problem changes less due to a smaller sampling variance in Algorithm~\ref{alg:varianceReduc}.
\corr{The variance for both approaches is depicted in Figure~\ref{fig:traceVariance} and listed in Table\ref{tab:Variances}. We want to emphasize that the approach based on solving the minimization problem is faster and simultaneously yields a lower variance in the entries of the self-diffusion compared to the estimation of the long-term limit.}

\begin{table}[!ht]
\begin{center}
\begin{tabular}{cccc}
\multicolumn{1}{l}{}    & \multicolumn{1}{l}{} & \multicolumn{1}{l}{$\umax{\ell \in \{0,\dots,N\}} \text{Var(Tr}(D_s({\ell/N})$))} & \multicolumn{1}{l}{$\frac{1}{N+1}\sum_{\ell = 0}^N \text{Var(Tr}(D_s({\ell/N})$))} \\ \hline
\multirow{2}{*}{N = 8}  & long-term limit      & $2.313 \cdot10^{-4}$                                                            & $1.072\cdot10^{-4}$                                          \\
                        & minimization         & $8.791\cdot10^{-9}$                                                               & $4.231\cdot10^{-9}$                                                         \\ \hline
\multirow{2}{*}{N = 15} & long-term limit      & $3.433\cdot10^{-4}$                                                            & $1.410\cdot10^{-4}$                                                            \\
                        & minimization         & $2.548\cdot10^{-4}$                                                               & $6.835\cdot10^{-5}$                                                          \\ \hline
\multirow{2}{*}{N=24}   & long-term limit      & $4.377\cdot10^{-4}$                                                               & $1.989\cdot10^{-4}$                                                           \\
                        & minimization         & $3.262\cdot10^{-4}$                                                               & $8.001\cdot10^{-5}$                                                           \\ \hline
\end{tabular}
\end{center}
\caption{\corr{We display the mean and maximum of the variance for the 100 approximations of the self-diffusion coefficient computed in Figure~\ref{fig:traceVariance}.
Note that we do not need to sample to evaluate $A^{\bm u}_{M,\ell}$ for  $N=8$. 
The variance for $N=8$ is solely caused by the random initializations of $R_{\bm s}^0$ in Algorithm\ref{alg:ALS}.
}
}
\label{tab:Variances}
\end{table}

\section{Conclusion}

\corr{Classical computational methods used for the approximation of the self-diffusion coefficient of a tagged particle process consist in exploiting the fact that the coefficient can be expressed as the long-time limit of the mean-square deviation of the tagged particle process~\eqref{eq:infiniteLimit}. These standard numerical methods consist in truncating the computational domain and approximate the long-time limit of the mean square deviation by an empirical average computed with a standard Monte Carlo approach from a large number of realizations of trajectories of the tagged particle. However, this classical approach suffers from two drawbacks: first, the amount of statistical noise due to the Monte-Carlo approach is very significant, as illustrated in our experiment, and yields larger errors than the error linked to the truncation of the computational domain. Indeed, it has been proved in~\cite{landim2002finite}, that the error in the approximation of the self-diffusion matrix decays exponentially with respect to the size of the finite computational domain. In contrast, the error due the statistical noise in Monte-Carlo averages only decays as the inverse of the square root of the number of random samples of trajectories. In addition, this traditional Monte Carlo approach requires to choose a priori a finite value of the time horizon, which induces additional errors which are difficult to estimate. }

In this work, we \corr{ have proposed a new approach to compute the self-diffusion matrix: we exploit the fact that the latter quantity can be expressed as the solution of a high-dimensional deterministic minimization problem~\eqref{eq:infinite} and we use tensor methods in order to build a numerical approximation of the solution of this problem. Here, the low-rank solutions are computed by means of an alternating optimization scheme.}

Our numerical experiments demonstrate that our approach is very advantageous compared to the classical Monte Carlo method. 
Indeed, with the same amount of computational resources, the variance in the obtained self-diffusion matrix is much smaller when using our new approach. \corr{Since the statistical noise is the leading source of errors in the approximation of this self-diffusion matrix, the low-rank approach proposed here provides a very interesting alternative to standard Monte-Carlo methods.}

\corr{This preliminary study yields several interesting remaining open questions, which we wish to investigate in future works. 
First, given the current numerical implementation of the low-rank algorithm proposed here, we are limited in terms of size of computational domains for which the present algorithm can be run because of round-off errors. 
We still believe that the approach described here is promising since, again, the errors linked to the truncation of the computational domain are small in comparison to statistical errors. 
However, we are not aware of any theoretical results on how to select the approximation rank, the number of ALS iterations and the number of samples depending on the size of the domain and the desired approximation error.
These parameters are crucial for the computational time required to approximate the self-diffusion matrix on larger domains.
In order to extend our approach to three-dimensional systems, it might be possible to combine the present method together with adapted domain decomposition approaches. We plan to study this in the future. 
}

\corr{Another very interesting question is the extension of the proposed approach to continuous-state diffusion processes, like Langevin dynamics in molecular dynamics for instance. We see the present work as a preliminary step towards the computation of diffusion coefficients out of Kinetic Monte-Carlo simulations, which can be fed with continuous space Molecular Dynamics simulations. 
An example of a possible track one could follow in this directions could be to introduce a discretization grid of the continuous space and compute approximate jump probability rates and occupancy probabilities of each cell of the discretization grid together with continuous state stochastic dynamics. This will require however, as a first step, the extension of the proposed approach to tagged particle systems where each grid site has the possibility to be occupied by more than one particle. We intend to investigate this issue in a future work. }

\subsection*{Acknowledgments}
\corr{This work was supported by the European Research Council under the European Union’s
Seventh Framework Programme (FP/2007-2013) / ERC Grant Agreement number 614492
and under the European Union’s Horizon 2020 Research and Innovation Programme, ERC
Grant Agreement number 810367, project EMC2.} The work was initiated during the CEMRACS 2021 summer school at CIRM, Luminy, Marseille. The authors also acknowledge funding by the ANR project COMODO (ANR-19-CE46-0002), the Center on Energy and Climate Change (E4C) and the I-Site FUTURE. 

The authors would like to thank Mi-Song Dupuy,  Guillaume Enchéry, Daniel Kressner and Dominik Schmid for stimulating discussions on this work.

{\footnotesize \input{main.bbl}}

\end{document}

%% file: 0Prae.tex
\usepackage{fullpage}
\usepackage[T1]{fontenc} 
\usepackage[latin1]{inputenc}
\usepackage[english]{babel}
\usepackage{selinput}
\SelectInputMappings{adieresis={ä},germandbls={ß}}

\usepackage{bm,bbm}
\usepackage{amsmath} 
\usepackage{amssymb} 
\usepackage{amsthm} 
\usepackage{enumerate}

\usepackage{graphicx} 
\usepackage{epstopdf}
\usepackage{subfig}
\usepackage[framed,numbered]{matlab-prettifier}
\usepackage{algorithm} 
\usepackage[noend]{algpseudocode}

\graphicspath{{./Figures/}}
\usepackage{multirow}
\usepackage{fancyvrb}
\usepackage{hyperref}
\usepackage{ mathrsfs }
\usepackage{tikz}
\usepackage[numbers,sort]{natbib}
\usetikzlibrary{shapes,arrows}

\newcommand{\R}{\mathbb{R}}
\newcommand{\bra}[1]{\left( #1 \right) }

\newcommand{\norm}[1]{\lVert #1 \rVert}
\newcommand{\abs}[1]{\left\lvert #1 \right\rvert}

\newcommand{\eps}{\varepsilon}

\DeclareMathOperator{\argmin}{argmin}

\newcommand{\umin}[1]{\underset{#1}{\min}\;}
\newcommand{\umax}[1]{\underset{#1}{\max}\;}

\newcommand{\corr}[1]{\textcolor{black}{#1}}

\newcommand{\matD}{\mathbb{D}}

\renewcommand{\bm}{\mathbf}

\newcommand{\bv}{\bm v}

\newcommand{\bu}{\bm u}
\newcommand{\by}{\bm y}
\newcommand{\bz}{\bm z}
\newcommand{\bs}{\bm s}
\newcommand{\bw}{\bm w}

\newcommand{\bEta}{\boldsymbol{\eta}}

\renewcommand{\bm}{\mathbf}

\theoremstyle{definition}
\newtheorem{rmrk}{Remark}[section]
